\documentclass{article}
\usepackage{amsfonts,amsmath,amssymb}
\usepackage{graphicx,epsfig,psfig}
\usepackage{checkend}
\usepackage{longtable,geometry}
\newtheorem{definition}{Definition}[section]
\begin{document}


\newtheorem{theorem}{Theorem}[section]
\newtheorem{proposition}{Proposition}[section]
\newtheorem{lemma}{Lemma}[section]
\newtheorem{corollary}{Corollary}[section]
\newtheorem{remark}{Remark}[section]
\newtheorem{proof}{Proof:}

\renewcommand{\thesection}{\arabic{section}}
\renewcommand{\theequation}{\thesection.\arabic{equation}}
\renewcommand{\thetheorem}{\thesection.\arabic{theorem}}
\numberwithin{equation}{section}
\numberwithin{theorem}{section}
\numberwithin{proposition}{section}
\numberwithin{lemma}{section}
\numberwithin{remark}{section}
\setcounter{secnumdepth}{5}


\newcommand{\cl}{\centerline}
\newcommand{\sms}{\smallskip}
\newcommand{\ms}{\medskip}
\newcommand{\bs}{\bigskip}
\newcommand{\noi}{\noindent}
\newcommand{\itl}[1]{\textit{#1}}
\newcommand{\blf}[1]{\textbf{#1}}
\newcommand{\dsty}{\displaystyle}
\newcommand{\txty}{\textstyle}
\newcommand{\ssty}{\scriptstyle}
\newcommand{\tty}{\texttt}


\newcommand\Par{\mathhexbox278\,}


\newcommand{\al}{\alpha}
\newcommand{\Al}{\Alpha}
\newcommand{\be}{\beta}
\newcommand{\Be}{\Beta}
\newcommand{\Gm}{\Gamma}
\newcommand{\gm}{\gamma}
\newcommand{\dl}{\delta}
\newcommand{\Dl}{\Delta}
\newcommand{\lm}{\lambda}
\newcommand{\Lm}{\Lambda}
\newcommand{\kp}{\kappa}
\newcommand{\varep}{\varepsilon}
\newcommand{\vp}{\varphi}
\newcommand{\sig}{\sigma}
\newcommand{\Sig}{\Sigma}
\newcommand{\om}{\omega}
\newcommand{\Om}{\Omega}
\newcommand{\uom}{\mbox{\boldmath$\omega$}}
\newcommand{\btau}{\mbox{\boldmath$\tau$}}
\newcommand{\bnu}{\mbox{\boldmath$\nu$}}
\newcommand{\up}{\upsilon}
\newcommand{\z}{\zeta}


\newcommand{\df}[1]{\buildrel\mbox{\small def}\over{#1}}
\newcommand{\op}[1]{\buildrel\mbox{\tiny o}\over{#1}}
\newcommand{\db}{\prime\prime}
\newcommand{\bsl}{\backslash}
\newcommand{\lb}{\lbrack\!\lbrack}
\newcommand{\rb}{\rbrack\!\rbrack}
\newcommand\la{\langle}
\newcommand\ra{\rangle}
\newcommand{\ev}{\equiv}
\newcommand{\nev}{\not\equiv}
\newcommand{\nn}{\mathbb{N}}
\newcommand{\qq}{\mathbb{Q}}
\newcommand{\zz}{\mathbb{Z}}
\newcommand{\rr}{\mathbb{R}}
\newcommand{\rn}{\rr^N}
\newcommand{\cc}{\mathbb{C}}
\newcommand{\id}{\mathbb{I}}
\newcommand{\bo}{\mathbb{O}}

\newcommand{\amsb}[1]{\mathbb{#1}}
\newcommand{\mcl}[1]{\mathcal{#1}}
\newcommand{\bl}[1]{\mathbf{#1}}
\newcommand{\ov}[1]{\overline{#1}}
\newcommand{\wt}[1]{\widetilde{#1}}
\newcommand{\wh}[1]{\widehat{#1}}

\newcommand{\lra}{\longrightarrow}
\newcommand{\LLR}{\Longleftrightarrow}
\newcommand{\LRA}{\Longrightarrow}
\newcommand{\LLA}{\Longleftarrow}


\newcommand{\bbox}{\vrule height.6em width.6em 
depth0em} 
\newcommand{\os}{\vbox{\hrule \hbox{\vrule 
height.6em depth0pt 
\hskip.6em \vrule height.6em depth0em}
\hrule}} 


\newcommand{\dvg}{\operatorname{div}}
\newcommand{\curl}{\operatorname{curl}}
\newcommand{\supp}{\operatorname{supp}}
\newcommand{\essup}{\operatornamewithlimits{ess\,sup}}
\newcommand{\essinf}{\operatornamewithlimits{ess\,inf}}
\newcommand{\essosc}{\operatornamewithlimits{ess\,osc}}
\newcommand{\osc}{\operatornamewithlimits{osc}}
\newcommand{\sign}{\operatorname{sign}}
\newcommand{\loc}{\operatorname{loc}}
\newcommand{\diam}{\operatorname{diam}}
\newcommand{\dist}{\operatorname{dist}}
\newcommand{\card}{\operatorname{card}}
\newcommand{\meas}{\operatorname{meas}}
\newcommand{\spn}{\operatorname{span}}
\newcommand{\dtm}{\operatorname{det}}
%


\newcommand{\overlim}{\mathop{\overline{\lim}}\limits}
\newcommand{\underlim}{\mathop{\underline{\lim}}\limits}
\newcommand{\ttop}[2]{\genfrac{}{}{0pt}{}{#1}{#2}}
\newcommand{\bcu}{\mathop{\txty{\bigcup}}\limits}
\newcommand{\bca}{\mathop{\txty{\bigcap}}\limits}
\newcommand{\bsu}{\mathop{\txty{\sum}}\limits}
\newcommand{\pro}{\mathop{\txty{\prod}}\limits}


\newcommand{\pl}{\partial}
\newcommand{\ptt}{\frac{\pl}{\pl t}}
\newcommand{\ppx}{\frac\pl{\pl x}}
\newcommand{\dds}{\frac d{ds}}
\newcommand{\ddt}{\frac d{dt}}


\newcommand{\intl}{\int\limits}
\newcommand{\iintl}{\iint\limits}


\def\Xint#1{\mathchoice
    {\XXint\displaystyle\textstyle{#1}}%
    {\XXint\textstyle\scriptstyle{#1}}%
    {\XXint\scriptstyle\scriptscriptstyle{#1}}%
    {\XXint\scriptscriptstyle\scriptscriptstyle{#1}}%
    \!\int}
\def\XXint#1#2#3{\setbox0=\hbox{$#1{#2#3}{\int}$}
    \vcenter{\hbox{$#2#3$}}\kern-0.5\wd0}
\def\bint{\Xint-}
\def\dashint{\Xint{\raise4pt\hbox to7pt{\hrulefill}}}


\newcommand{\ovl}[3]{\int_{#1}^{#2}\kern-#3pt\raise4pt\hbox to7pt{\hrulefill}\ }

\newcommand{\ovll}[3]{\intl_{#1}^{#2}\kern-#3pt\raise4pt\hbox to7pt{\hrulefill}\ }

\newcommand{\tvl}[2]{\iint_{#1}\kern-#2pt\raise4pt\hbox to7pt{\hrulefill}\ }



\newcommand{\omt}{\Om_T}
\newcommand{\plo}{\partial\Omega}
\newcommand{\ovo}{\bar{\Om} }

%
\newcommand{\ci}[1]{C^\infty\!\left({#1}\right)}
\newcommand{\cio}[1]{C_o^\infty\!\left({#1}\right)}
\newcommand{\lloc}[1]{L_{\loc}\!\left({#1}\right)}
\newcommand{\xy}{|x-y|}


\newcommand{\intom}{\intl_{\Om}}
\newcommand{\intbo}{\intl_{\plo}}
\newcommand{\inom}{\int_{\Om}}
\newcommand{\inbo}{\int_{\plo}}
\newcommand{\intrn}{\intl_{\rn}}


\newcommand{\bye}{\end{document}}



%
%
%

%
%
%
%
%

%
%

\input harnack.mac

\begin{center}  {\huge\textbf{On the Optimal Control of the Free Boundary Problems for the 
Second Order Parabolic Equations. I.Well-posedness and Convergence of the Method of Lines }}
\par \medskip\bigskip\end{center}
\begin{center} {\Large\textsc{Ugur G. Abdulla}}
\par \medskip\bigskip\end{center}
\begin{center} {\large\noindent \textsc{Department of Mathematics, Florida Institute of Technology, Melbourne, Florida 32901}}
\par \medskip\bigskip\end{center}
{\bf Abstract.} We develop a new variational formulation of the inverse Stefan problem, where information on the heat flux on the fixed boundary is missing and must be found along with the temperature and free boundary. We employ optimal control framework, where boundary heat flux and free boundary are components of the control vector, and optimality criteria consists of the minimization of the sum of $L_2$-norm declinations from the available measurement of the temperature flux on the fixed boundary and available 
information on the phase transition temperature on the free boundary. This approach allows 
one to tackle situations when the phase transition temperature is not known explicitly, and is available through measurement with possible error. It also allows for the development of iterative numerical methods of least computational cost due to the fact that for every given control vector, the parabolic PDE is solved in a fixed region instead of full free boundary problem. We prove well-posedness in Sobolev spaces framework and 
convergence of discrete optimal control problems to the original problem both with respect to cost functional and control.



\newpage
\section{Description of Main Results}\label{description of results,historical remarks}
\subsection{Introduction and Motivation}\label{E:1:1}
{\large
Consider the general one-phase Stefan problem (\cite{Friedman1, Meyrmanov}): find the temperature function $u(x,t)$ and the free boundary $x=s(t)$ from the following conditions
\begin{equation}\label{Eq:W:1:1}
(a(x,t)u_{x})_{x}+b(x,t)u_{x}+c(x,t)u-u_{t}=f(x,t),\quad \text{for}~(x,t) \in \Omega
\end{equation}
\begin{equation}\label{Eq:W:1:2}
u(x,0)=\phi(x),\qquad 0 \leq x \leq s(0)=s_{0}
\end{equation}
\begin{equation}\label{Eq:W:1:3}
a(0,t)u_{x}(0,t)=g(t),\qquad 0 \leq t \leq T
\end{equation}
\begin{equation}\label{Eq:W:1:4}
a(s(t),t)u_{x}(s(t),t) + \gamma(s(t),t)s'(t)=\chi(s(t),t),\qquad 0 \leq t \leq T
\end{equation}
\begin{equation}\label{Eq:W:1:5}
u(s(t),t)=\mu(t),\qquad 0 \leq t \leq T
\end{equation}
where $a$, $b$, $c$, $f$, $\phi$, $g$, $\gamma$, $\chi$, $\mu$ are known functions and
\begin{equation}\label{Eq:W:1:6}
a(x,t)\geq a_{0}>0,\quad s_{0}>0
\end{equation}
\[ \Omega=\left\{ (x,t) :~0<x<s(t),~0<t\leq T \right\} \]
In the physical context, $f$ characterizes the density of the sources, $\phi$ is the initial temperature, $g$ is the heat flux on the fixed boundary and $\mu$ is the phase transition temperature.  

Assume now that some of the data is not available, or involves some measurement error.  For example, assume that the heat flux $g(t)$ on the fixed boundary $x=0$ is not known and must be found along with the temperature $u(x,t)$ and the free boundary $s(t)$.  In order to do that, some additional information is needed.  Assume that this additional information is given in the form of the temperature measurement along the boundary $x=0$:
\begin{equation}\label{Eq:W:1:7}
u(0,t)=\nu(t),\quad \text{for}~0 \leq t \leq T
\end{equation}
\textbf{Inverse Stefan Problem (ISP):} \textit{Find the functions $u(x,t)$ and $s(t)$ and the boundary heat flux $g(t)$ satisfying conditions (\ref{Eq:W:1:1})-(\ref{Eq:W:1:7}}).

ISP is not well posed in the sense of Hadamard.  If there is no coordination between the input data, the exact solution may not exist.  Even if it exists, it might be not unique, and most importantly, there is no continuous dependence of the solution on the data. 
Inverse Stefan problem was first mentioned in \cite{Cannon3}, in the form of finding
a heat flux on the fixed boundary which provides a desired free boundary. This problem is similar to non-characteristic Cauchy problem for the heat equation. The variational approach
for solving this ill-posed inverse Stefan problem was performed in \cite{BudakVasileva1, BudakVasileva2}. First result on the optimal control of the Stefan problem appeared in \cite{Vasilev}. It consists of finding optimal value of the external temperature along the fixed boundary, in order to ensure that the solutions of the Stefan problem are close to the measurements taken
at the final moment. In \cite{Vasilev} existence result was proved. In \cite{Yurii} the Frechet differentiability and the convergence of the difference schemes was proved for the same problem and Tikhonov regularization was suggested. Later development of the inverse Stefan problem was along these two lines: Inverse Stefan problems with given phase boundaries were considered in \cite{Alifanov,Bell,Budak,Cannon,Carasso,Ewing1,Ewing2,Hoffman,Sherman,Goldman}; optimal control of Stefan problems, or equivalently inverse problems with unknown phase boundaries were investigated in \cite{Baumeister,Fasano,Hoffman1,Hoffman2,Jochum2,Jochum1,Knabner,Lurye,Nochetto, Niezgodka,Primicero,Sagues,Talenti,Goldman}. We refer to monography \cite{Goldman} for a complete list of references of both types of inverse Stefan problems, both for linear and quasilinear parabolic equations. 
The main methods used to solve inverse Stefan problem are based on variational formulation, method of quasi-solutions or Tikhonov regularization which takes into account ill-posedness in terms of the dependence of the solution on the inaccuracy involved in the measurement (\ref{Eq:W:1:7}), Frechet differentiability and iterative conjugate gradient methods for numerical solution. Despite its effectiveness, this approach has some deficiencies in many practical applications:
\begin{itemize}
\item Solution of the inverse Stefan problem is not continuously dependent on the phase transition temperature $\mu(t)$: small perturbation of the phase transition temperature may imply significant change of the solution to the inverse Stefan problem. Accordingly, any regularization which equally takes into account instability with respect to both $\nu(t)$ from measurement (\ref{Eq:W:1:7}), and the phase
transition temperature $\mu(t)$ from (\ref{Eq:W:1:5}) will be preferred. It should be also mentioned that in many applications the phase transition temperature is not known explicitly. In many processes the melting temperature of pure material at a given external action depends on the process evolution. For example, gallium (Ga, atomic number 31) may remain in the liquid phase at temperatures well below its mean melting temperature (\cite{Meyrmanov}). 
\item Numerical implementation of the iterative gradient type methods within the existing approach requires to solve full free boundary problem at every step of 
the iteration, and accordingly requires quite high computational cost. Iterative gradient method which requires at every step solution of the boundary value problem in a fixed region would definitely be much more effective in terms of the computational cost.  
\end{itemize}

The main goal of this project is to develop a new variational approach based on the optimal control theory which is capable of addressing both of the mentioned issues and allows the inverse Stefan problem to be solved numerically with least computational cost by using conjugate gradient methods in Hilbert spaces. 
In this paper we prove the existence of the optimal control and convergence of the family of time-discretized optimal control problems to the continuous problem both with respect to cost functional and control. We employ Sobolev spaces framework which allows to reduce the reguarity and structural requirements on the data. We address the problems of convergence of the fully discretized family
of optimal control problems, Frechet differentiability and iterative conjugate gradient methods in Hilbert spaces in an upcoming paper.

Throughout the paper we use usual notation for Sobolev spaces according to references \cite{LSU,BIN,Nikolski,Solonnikov1,Solonnikov2}.

In the next section we formulate a new variational formulation of the inverse problem which 
takes into account the described deficiencies.

\subsection{Optimal Control Problem}\label{E:1:2}
Consider a minimization of the cost functional
\begin{equation}
\mathcal{J}(v)=\beta_{0}\Vert u(0,t)-\nu(t)\Vert_{L_{2}[0,T]}^{2}+\beta_{1}\Vert u(s(t),t)-\mu(t)\Vert_{L_{2}[0,T]}^2\label{Eq:W:1:8}
\end{equation}
on the control set 
\begin{equation*}
V_{R}=\left\{ v=(s,g) \in W_{2}^{2}[0,T]\times W_{2}^{1}[0,T]: \delta \leq s(t)\leq l, s(0)=s_0,
\max (~\Vert s\Vert_{W_{2}^{2}}; ~\Vert g\Vert_{W_{2}^{1}} \leq R\right\}
\end{equation*}
where $\delta, l,R, \beta_0, \beta_1$ are given positive numbers, and $u=u(x,t;v)$ be a solution of the Neumann problem (\ref{Eq:W:1:1})-(\ref{Eq:W:1:4}).   

\begin{definition}
The function $u \in W_{2}^{1,1}(\Omega)$ is called a weak solution of the problem (\ref{Eq:W:1:1})-(\ref{Eq:W:1:4}) if $u(x,0)=\phi(x) \in W_{2}^{1}[0,s_0]$ and
\begin{gather}
0=\int_{0}^{T}\int_{0}^{s(t)}[ a u_{x}\Phi_{x}-bu_{x}\Phi - c u \Phi + u_{t} \Phi+f\Phi] \,dx\,dt \nonumber\\
 +\int_{0}^{T}[ \gamma(s(t),t)s'(t)-\chi(s(t),t)]\Phi(s(t),t)\, dt
+\int_{0}^{T}g(t)\Phi(0,t)\, dt\label{Eq:W:1:9}
\end{gather}
for arbitrary $\Phi \in W_{2}^{1,1}(\Omega)$
\end{definition}
We also need a notion of weak solution from $V_{2}(\Omega)$ of the Neumann problem:
\begin{definition}
The function $u \in V_{2}(\Omega)$ is called a weak solution of (\ref{Eq:W:1:1})-(\ref{Eq:W:1:4}) if
\begin{gather}
0=\int_{0}^{T}\int_{0}^{s(t)}[ a u_{x}\Phi_{x}-bu_{x}\Phi - c u \Phi - u \Phi_{t}+f\Phi] \,dx\,dt -\int_{0}^{s_0}\phi(x)\Phi(x,0)\,dx+\nonumber\\
 \int_{0}^{T}g(t)\Phi(0,t)\, dt+\int_{0}^{T}[ \gamma(s(t),t)s'(t)-u(s(t),t)s'(t)-\chi(s(t),t)]\Phi(s(t),t)\, dt \label{Eq:W:1:10}
\end{gather}
for arbitrary $\Phi \in W_{2}^{1,1}(\Omega)$ such that $\left. \Phi\right|_{t=T}=0$.
\end{definition}
If $u$ is a weak solution either from $V_2(\Omega)$ (or $W_{2}^{1,1}(\Omega)$), then traces $\left. u\right|_{x=0}$ and
$\left. u\right|_{x=s(t)}$ are elements of $L_2[0,T]$, when $s\in W_2^2[0,T]$ (\cite{Nikolski, LSU}) and cost functional $\mathcal{J}(v)$
is well defined. Furthermore, formulated optimal control problem will be called Problem $I$.
\subsection{Discrete Optimal Control Problem}\label{E:1:3}
Let
\[\omega_{\tau}=\{ t_{j}=j \cdot \tau,~j=0,1,\ldots,n\} \]
be a grid on $[0,T]$ and $\tau=\frac{T}{n}$. Consider a discretized control set
\begin{equation*}
V^n_{R}=\{ [v]_{n}=([s]_n,[g]_n) \in {\mathbb R}^{2n+2}:~0<\delta\leq s_{k} \leq l,~ \max(\Vert [s]_{n}\Vert_{w_{2}^{2}}^2; ~\Vert [g]_{n}\Vert_{w_{2}^{1}}^2) \leq R^2\}
\end{equation*}
where,
\[ [s]_n=(s_0,s_1,...,s_n) \in {\mathbb R}^{n+1}, \ [g]_n=(g_0,g_1,...,g_n) \in {\mathbb R}^{n+1} \]
\[ 
\Vert [s]_{n}\Vert_{w_{2}^{2}}^2= \sum\limits_{k=0}^{n-1}\tau s_k^2+\sum\limits_{k=1}^{n}\tau s_{\overline{t},k}^2+\sum\limits_{k=1}^{n-1}\tau s_{\overline{t}t,k}^2, \  \Vert [g]_{n}\Vert_{w_{2}^{1}}^2= \sum\limits_{k=0}^{n-1}\tau g_k^2+\sum\limits_{k=1}^{n}\tau g_{\overline{t},k}^2.
\]
under the standard notation for the finite differences:
\[ s_{\overline{t},k}=\frac{s_k-s_{k-1}}{\tau}, \ s_{t,k}=\frac{s_{k+1}-s_{k}}{\tau}, \ s_{\overline{t}t,k}^2=\frac{s_{k+1}-2s_k+s_{k-1}}{\tau^2}. \]
Introduce two mappings $\mathcal{Q}_n$ and $\mathcal{P}_n$ between continuous and discrete control sets:
\[ \mathcal{Q}_n(v)=[v]_n=([s]_n,[g]_n), \quad \text{for}~ v\in V_R \]
where $s_k=s(t_k), g_k=g(t_k), k=0,1,...,n$.
\[ \mathcal{P}_n([v]_n)=v^n=(s^n,g^n)\in W_2^2[0,T]\times W_2^1[0,T] \quad \text{for}~ [v]_n \in V_R^n, \]
where
\begin{equation}\label{Eq:W:1:11}
s^n(t)=
\left\{
\begin{array}{l}
s_0+\frac{t^2}{2\tau} s_{\overline{t},1} \ \ 0\le t \le \tau,\\
s_{k-1}+(t-t_{k-1}-\frac{\tau}{2})s_{\overline{t},k-1}+\frac{1}{2}(t-t_{k-1})^2 s_{\overline{t}t,k-1} \ \ t_{k-1}\le t \le t_k, k=\overline{2,n}.
\end{array}\right.
\end{equation}
\begin{equation*}
g^n(t)=g_{k-1}+\frac{g_k-g_{k-1}}{\tau}(t-t_{k-1}), \ \ t_{k-1} \le t \le t_k, k=\overline{1,n}.
\end{equation*}
Introduce Steklov averages
\[ d_{k}(x)=\frac{1}{\tau}\int_{t_{k-1}}^{t_{k}}d(x,t)\,dt, \ h_{k}=\frac{1}{\tau}\int_{t_{k-1}}^{t_{k}}h(t)\,dt, \]
where $d$ stands for any of the functions $a$, $b$, $c$, $f$, and $h$ stands for any of the functions $\nu$, $\mu$. Given $v=(s,g) \in V_R$ we define Steklov averages of traces
\begin{equation}\label{Eq:W:1:12}
\chi^{k}_s=\frac{1}{\tau} \int_{t_{k-1}}^{t_{k}}\chi(s(t),t) \,dt, \
(\gamma_s s')^k=\frac{1}{\tau} \int_{t_{k-1}}^{t_{k}}\gamma(s(t),t)s'(t) \,dt.
\end{equation}
Given $[v]_n=([s]_n,[g]_n) \in V_R^n$ we define Steklov averages $\chi^{k}_{s^n}$ and $(\gamma_{s^n} (s^n)')^k$ through (\ref{Eq:W:1:12}) with  $s$ replaced by $s^n$ from (\ref{Eq:W:1:11}).

Next we define a discrete state vector through time-discretization of the integral identity (\ref{Eq:W:1:9})
\begin{definition}\label{discretestatevector}
Given discrete control vector $[v]_n$, the vector function
\[ [u([v]_n)]_n=(u(x;0),u(x;1),...,u(x;n))  \]
is called a discrete state vector if 
\begin{description}
\item{\bf(a)} $u(x;0)=\phi(x)\in W_2^1[0,s_0]$; 
\item{\bf(b)} For arbitrary $k=1,2,\ldots,n$, $u(x;k) \in W_{2}^{1}[0,s_{k}]$ satisfy the integral identity
\begin{align}
&\int_{0}^{s_{k}}\Big ( a_{k}(x)\frac{d u(x;k)}{d x} \eta'(x)-b_{k}\frac{d u(x;k)}{d x} \eta(x) - c_{k}(x)u(x;k) \eta(x) + f_{k}(x)\eta(x)  \nonumber\\
&\qquad + u_{\overline{t}}(x;k) \eta(x) \Big ) \,dx + \Big ( (\gamma_{s^n} (s^n)')^k-\chi^{k}_{s^n} \Big )\eta(s_{k})+g_{k}\eta(0)=0,\label{Eq:W:1:13}
\end{align}
for arbitrary $\eta \in W_{2}^{1}[0,s_{k}]$, where
\[ u_{\overline{t}}(x;k)=\frac{u(x;k)-u(x;k-1)}{\tau}. \]
\item{\bf(c)} For arbitrary $k=0,1,...,n$, $u(x;k)\in W_2^1[0,s_k]$ iteratively continued to $[0,l]$ as
\begin{equation}\label{Eq:W:1:14}
u(x;k)=u(2^ns_k-x;k), \ 2^{n-1}s_k \le x \le 2^ns_k, n=\overline{1,n_k}, \ n_k\le N=1+\log_2\Big [ \frac{l}{\delta}\Big ]
\end{equation}
where $[r]$ means integer part of the real number $r$.
\end{description}
\end{definition}

Consider a discrete optimal control problem of minimization of the cost functional
\begin{equation}\label{Eq:W:1:15}
\mathcal{I}_n([v]_n)=\beta_{0}\tau\sum\limits_{k=1}^n \Big ( u(0;k)-\nu_k \Big )^2+\beta_{1}\tau\sum\limits_{k=1}^n \Big ( u(s_k;k)-\mu_k \Big )^2
\end{equation}
on a set $V_R^n$ subject to the state vector defined in Definition 1.3. Furthermore, formulated discrete optimal control problem will be called Problem $I_n$.

Throughout we use piecewise constant and piecewise linear interpolations of the discrete state vector:
given discrete state vector $[u([v]_n)]_n=(u(x;0),u(x;1),...,u(x;n))$, let
\[ u^\tau(x,t)=u(x;k), \quad \text{if}~ t_{k-1}<t\le t_k, \ 0\le x \le l, \ k=\overline{0,n}, \]
\[ \hat{u}^\tau(x,t)=u(x;k-1)+u_{\overline{t}}(x;k)(t-t_{k-1}), \quad \text{if}~ t_{k-1}<t\le t_k, \ 0\le x \le l, \ k=\overline{1,n}, \]
\[ \hat{u}^\tau(x,t)= u(x;n), \quad \text{if}~ t\ge T, \ 0\le x \le l. \]
Obviously, we have
\[ u^\tau \in V_2^{1,0}(D), \ \ \hat{u}^\tau \in W_2^{1,1}(D)  \]

\subsection{Formulation of the Main Result}\label{E:1:4}
Let 
\[ D=\{ (x,t):~0<x<l,~0<t\leq T\}\]
Throughout the whole paper, with the exeption of Section~\ref{firstenergyestimate}, we assume the following conditions are satisfied by the data:
\[ a,b,c \in L_{\infty}(D), \ f \in L_2(D), \]
\[ \phi \in W_2^1[0,s_0], \ \gamma, \chi \in W_2^{1,1}(D), \ \mu,\nu \in L_2[0,T], \]
the coefficient $a$ satisfies (\ref{Eq:W:1:6}) almost everywhere on $D$, 
the generalized derivatives $\frac{\partial a}{\partial t},  \frac{\partial a}{\partial x}$ exists and 
\begin{equation}\label{conditionon_a}
\frac{\partial a}{\partial x}\in L_{\infty}(D), \ \ \int_{0}^{T}esssup_{0 \leq x \leq l} \left| \frac{\partial a}{\partial t}\right| dt <+\infty.
\end{equation}
Our main theorems read:
\begin{theorem}\label{existence}
The Problem $I$ has a solution, i.e.
\[ V_*=\{v\in V_R: \mathcal{J}(v)= \mathcal{J}_* \equiv \inf\limits_{v\in V_R} \mathcal{J}(v) \} \neq \emptyset \]
\end{theorem}
\begin{theorem}\label{convergence}
Sequence of discrete optimal control problems $I_n$ approximates the optimal control problem $I$ with respect to functional, i.e.
\begin{equation}\label{Eq:W:1:15}
\lim\limits_{n\to +\infty} \mathcal{I}_{n_*}=\mathcal{J}_*, 
\end{equation}
where
\[ \mathcal{I}_{n_*}=\inf\limits_{V_R^n} \mathcal{I}_n([v]_n), \ n=1,2,... \]
If $[v]_{n_\ep}\in V_R^n$ is chosen such that
\[ \mathcal{I}_{n_*} \le \mathcal{I}_n([v]_{n_\ep})\le \mathcal{I}_{n_*}+\ep_n, \ \ep_n \downarrow 0, \]
then the sequence $v_n=(s_n,g_n)=\mathcal{P}_n([v]_{n_\ep})$ converges to some element $v_*=(s_*,g_*) \in V_*$ weakly
in $W_2^2[0,T] \times W_2^1[0,T]$, and strongly in $W_2^1[0,T] \times L_2[0,T]$. In particular $s_n$ converges to $s_*$ uniformly on $[0,T]$. Moreover, piecewise linear interpolation $\hat{u}^\tau$ of the 
discrete state vector $[u[v]_{n_\ep}]_n$ converges to the solution $u(x,t;v_*) \in W_2^{1,1}(\Omega_*)$ of the Neumann 
problem (\ref{Eq:W:1:1})-(\ref{Eq:W:1:4}) weakly in $W_2^{1,1}(\Omega_*)$.
\end{theorem} 
\section{Preliminary Results}\label{preliminaries}
In a Lemma~\ref{existencediscretestate} below we prove existence and uniqueness 
of the discrete state vector $[u([v]_n)]_n$ (see Definition~\ref{discretestatevector}) for arbitrary discrete control vector $[v]_n \in V_R^n$.
In a Lemma~\ref{generalcriteria} we remind a general approximation criteria for the optimal control problems from   (\cite{Vasilev1}). In a Lemma~\ref{mappings} we prove some properties of the mappings
$\mathcal{Q}_n$ and $\mathcal{P}_n$ between continuous and discrete control sets.

\begin{lemma}\label{existencediscretestate}
For sufficiently small time step $\tau$, there exists a unique discrete state vector $[u([v]_n)]_n$ for arbitrary discrete control vector $[v]_n \in V_R^n$.
\end{lemma}

Proof. To prove uniqueness, it is enough to show that if 
\[ u(x;k-1) \equiv 0, \ (\gamma_{s^n} (s^n)')^k =0, \ \chi^{k}_{s^n}=0, \ g^k=0, \ f_k(x) \equiv 0 \]
then $u(x;k)$ which solves (\ref{Eq:W:1:13}) vanishes identically. Under these assumptions by choosing $\eta(x)=u(x;k)$ in (\ref{Eq:W:1:13}) we have
\begin{equation}\label{Eq:W:2:1}
\int_{0}^{s_{k}}\Big ( a_{k}(x)\Big ( \frac{d u(x;k)}{d x}\Big )^2 -b_{k}\frac{d u(x;k)}{d x} u(x;k) - c_{k}(x)u^2(x;k)+ \frac{1}{\tau} u^2(x;k)\Big ) \,dx =0.
\end{equation}
Using (\ref{Eq:W:1:6}) and Cauchy inequality with $\ep >0$ we derive that
\begin{gather}
a_0\int_{0}^{s_{k}}\Big ( \frac{d u(x;k)}{d x}\Big )^2 \,dx + \frac{1}{\tau}\int_{0}^{s_{k}}u^2(x;k) \,dx \le  \nonumber\\
 \frac{\ep M}{2} \int_{0}^{s_{k}}\Big ( \frac{d u(x;k)}{d x}\Big )^2 \,dx + \Big ( \frac{M}{2\ep}+M \Big ) \int_{0}^{s_{k}}u^2(x;k) \,dx,\label{Eq:W:2:2}
\end{gather}
where
\[ M=\max \Big ( ||a||_{L_{\infty}(D)}; ||b||_{L_{\infty}(D)}; ||c||_{L_{\infty}(D)} \Big ). \]
By choosing $\ep=a_0/M$ in (\ref{Eq:W:2:2}) we have
\begin{equation}\label{Eq:W:2:3}
\frac{a_0}{2}\int_{0}^{s_{k}}\Big ( \frac{d u(x;k)}{d x}\Big )^2 \,dx + \Big ( \frac{1}{\tau}-\frac{1}{\tau_0} \Big )\int_{0}^{s_{k}}u^2(x;k) \,dx \le 0, 
\end{equation}
where
\[ \tau_0 = \Big ( \frac{M^2}{2a_0}+M \Big )^{-1}. \]
From (\ref{Eq:W:2:3}) it follows that $u(x;k)\equiv 0$ if $\tau < \tau_0$.

To prove an existence we apply Galerkin method. Consider an approximate solution
\[u_{N}(x)=\sum_{i=1}^{N}d_{i}\psi_{i}(x)\]
where $\{\psi_{i}\}$ is a fundamental system in $W_{2}^{1}[0,s_{k}]$ and the coefficients $\{d_{i}\}$ solve the following system
\begin{gather}
\int_{0}^{s_{k}}\Big[a_{k}(x)\frac{d u_{N}}{d x} \psi'_{i}(x) - b_{k}(x)\frac{d u_{N}}{d x}\psi_{i}(x)- c_{k}(x)u_{N}(x) \psi_{i}(x)+\frac{1}{\tau}u_{N}(x)\psi_{i}(x)
 + f_{k}(x) \psi_{i}(x) \Big]\, dx\nonumber\\
=\frac{1}{\tau}\int_{0}^{s_{k}}u(x;k-1)\psi_{i}(x)\, dx -
\left[ (\gamma_{s^n} (s^n)')^k-\chi^{k}_{s^n} \right] \psi_{i}(s_{k})
 - g^{k}\psi_{i}(0) \ i=1,\ldots,N\label{eq:galerkin-sys}
\end{gather}
which is equivalent to
\begin{gather}
\sum_{j=1}^{N}\int_{0}^{s_{k}}\Big [ a_{k}(x)\psi'_{j}(x) \psi'_{i}(x)- b_{k}(x)\psi_{j}'(x)\psi_{i}(x)- c_{k}(x)\psi_{j}(x) \psi_{i}(x)+\frac{1}{\tau}\psi_{j}(x) \psi_{i}(x)\Big ]\, dx \ d_j= \nonumber\\
 \int_{0}^{s_{k}}\Big [ -f_{k}(x) \psi_{i}(x)+ \frac{1}{\tau}u(x;k-1)\psi_{i}(x)\Big ]-\left[ (\gamma_{s^n} (s^n)')^k-\chi^{k}_{s^n} \right] \psi_{i}(s_{k})
 - g^{k}\psi_{i}(0) \label{eq:galerkin-sys-equiv}
\end{gather}
$i=1,\ldots,N.$ Homogeneous system corresponding to
\eqref{eq:galerkin-sys-equiv} is 
\begin{gather}
\sum_{j=1}^{N}\int_{0}^{s_{k}}\Big [ a_{k}(x)\psi'_{j}(x) \psi'_{i}(x)- b_{k}(x)\psi_{j}'(x)\psi_{i}(x)- c_{k}(x)\psi_{j}(x) \psi_{i}(x)\nonumber\\ +\frac{1}{\tau}\psi_{j}(x) \psi_{i}(x)\Big ]\, dx \ d_j=0, \ 
\ i=1,2,\ldots,N\label{eq:galerkin-sys-homo}
\end{gather}
Let us multiply each equation in \eqref{eq:galerkin-sys-homo} by $d_{i}$ and add with respect to $i$:
\begin{equation}
\int_{0}^{s_{k}}\Big [ a_{k}(x)\left(\frac{d u_{N}(x)}{d x} \right)^{2}- b_{k}(x)\frac{d u_{N}(x)}{d x}u_{N}(x)dx-c_{k}(x)u_{N}^2(x) +\frac{1}{\tau}u_{N}^2(x))\Big ]\, dx =0\label{eq:galerkin-sys-homo-sum}
\end{equation}
As before, from \eqref{eq:galerkin-sys-homo-sum} it follows that $u_{N}\equiv 0$, and therefore the homogeneous system \eqref{eq:galerkin-sys-homo} has only the trivial solution.  This proves the uniqueness of the approximate solution $u_{N}(x)$.  Let us now prove uniform estimation of the sequence $\{u_{N}(x)\}$.  Multiply \eqref{eq:galerkin-sys} by $d_{i}$ and add with respect to $i=1,\ldots,N$:
\begin{gather}
\int_{0}^{s_{k}}\Big[a_{k}(x)\Big (\frac{d u_{N}}{d x}\Big )^2 - b_{k}(x)\frac{d u_{N}}{d x}u_N(x)- c_{k}(x)u_{N}^2(x)+\frac{1}{\tau}u_{N}^2(x)
 + f_{k}(x) u_N(x) \Big]\, dx\nonumber\\
=\frac{1}{\tau}\int_{0}^{s_{k}}u(x;k-1)u_N(x)\, dx -
\left[ (\gamma_{s^n} (s^n)')^k-\chi^{k}_{s^n} \right] u_N(s_{k})
 + g^{k}u_N(0).\label{eq:galerkin-sys2}
\end{gather}
We estimate the four integrals on the left-hand side of \eqref{eq:galerkin-sys2}
as we did before to prove (\ref{Eq:W:2:3}) and derive 
\begin{gather}
\frac{a_0}{2}\int_{0}^{s_{k}}\Big ( \frac{d u_N(x)}{d x}\Big )^2 \,dx +  \frac{1}{2\tau}\int_{0}^{s_{k}}u_N^2(x) \,dx \le |g^{k}||u_{N}(0)|+ \nonumber\\ \left[ |(\gamma_{s^n} (s^n)')^k|+|\chi^{k}_{s^n}|\right]|u_{N}(s_{k})|+
\int_{0}^{s_{k}}\left[ |f_{k}(x)|+
\frac{1}{\tau}|u(x;k-1)|\right ]|u_{N}(x)| \,dx \label{eq:galerkin-sys-est}
\end{gather}
for all $\tau \le \frac{\tau_0}{2}$.
By Morrey's inequality we have
\begin{equation}
\max\{|u_N(0)|; |u_N(s_k)|\}\le \Vert u_N \Vert_{C[0,s_k]}\le C \Vert u_N \Vert_{W_2^1[0,s_k]},\label{Morrey}
\end{equation}
where the constant $C$ is independent of $N$ and $\tau$. By using Cauchy inequalities with appropriately chosen $\ep >0$, from \eqref{eq:galerkin-sys-est} and \eqref{Morrey} 
it easily follows that
\begin{align}
\Vert u_N \Vert_{W_{2}^{1}[0,s_{k}]}^{2} &\leq C \Big ( \Vert u(x;k-1)\Vert_{L_{2}[0,s_{k}]} +\Vert f_{k}\Vert_{L_{2}(0,s_{k})} +\big|(\gamma_{s^n} (s^n)')^{k}\big|^{2} 
 +|\chi_{s^n}^{k}|^{2}+|g^{k}|^{2}\Big ) \label{eq:lemma1-c3}
\end{align}
where $C$ does not depend on $N$, but depends on the time step $\tau$.
From \eqref{eq:lemma1-c3} it follows that $\{ u_{N}\}$ is weakly compact in $W_{2}^{1}[0,s_{k}]$.  Let $v(x)$ be its weak limit point in $W_{2}^{1}[0,s_{k}]$.  Passing to the limit in \eqref{eq:galerkin-sys} it follows that $v(x)$ satisfies (\ref{Eq:W:1:13}) for $\eta(x)=\psi_{i}(x)$.  Since $\{ \psi_{i}\}$ is a fundamental system in $W_{2}^{1}[0,s_{k}]$, it follows that $v(x)$ satisfies (\ref{Eq:W:1:13}) for every $\eta(x) \in W_{2}^{1}[0,s_{k}]$.  Hence $v(x)=u(x;k)$ is a solution of (\ref{Eq:W:1:13}) and in view of uniqueness the whole sequence $u_{N}$ converges weakly in $W_{2}^{1}[0,s_{k}]$ to $u(x;k)$. Lemma is proved.

The following known criteria will be used in the proof of Theorem~\ref{convergence}.  

\begin{lemma}\label{generalcriteria}
\cite{Vasilev1} Sequence of discrete optimal control problems $I_n$ approximates the
continuous optimal control problem $I$ if and only if the following conditions
are satisfied:
\begin{description}
\item{\bf(1)} for arbitrary sufficiently small $\ep>0$ there exists number $N_1=N_1(\ep)$ such that $\mathcal{Q}_N(v)\in V^n_R$ for all $v \in V_{R-\ep}$ and $N\ge N_1$; and for any fixed $\ep>0$ and for all $v\in V_{R-\ep}$ the following inequality is satisfied:
\begin{equation}\label{firstcondition}
\limsup\limits_{N\to \infty} \Big ( \mathcal{I}_N(\mathcal{Q}_N(v))-\mathcal{J}(v) \Big ) \le 0.
\end{equation}
\item{\bf(2)} for arbitrary sufficiently small $\ep>0$ there exists number $N_2=N_2(\ep)$ such that $\mathcal{P}_N([v]_N)\in V_{R+\ep}$ for all $[v]_N \in V^N_R$ and $N\ge N_2$; and for all $[v]_N\in V^N_R$, $N\ge 1$ the following inequality is satisfied:
\begin{equation}\label{secondcondition}
\limsup\limits_{N\to \infty} \Big ( \mathcal{J}(\mathcal{P}_N([v]_N))-\mathcal{I}_N([v]_N)  \Big ) \le 0.
\end{equation}
\item{\bf(3)} the following inequalities are satisfied:
\begin{align}
\limsup\limits_{\ep \to 0} \mathcal{J}_*(\ep) \ge \mathcal{J}_*, 
\ \ \liminf\limits_{\ep \to 0} \mathcal{J}_*(-\ep) \le \mathcal{J}_*,
\end{align}\label{thirdcondition}
where $\mathcal{J}_*(\pm\ep)=\inf\limits_{V_{R\pm \ep}}\mathcal{J}(u)$.
\end{description}
\end{lemma}
In the next lemma we prove that the mappings $\mathcal{Q}_n$ and $\mathcal{P}_n$ introduced in Section~\ref{E:1:3} satisfy the conditions of Lemma~\ref{generalcriteria}.
\begin{lemma}\label{mappings}
For arbitrary sufficiently small $\ep > 0$ there exists $n_{\ep}$ such that
\begin{equation}\label{Eq:W:2:15}
\mathcal{Q}_n(v) \in V_R^n, \quad \text{for all}~v\in V_{R-\ep} \quad \text{and}~ n>n_\ep.
\end{equation}
\begin{equation}\label{Eq:W:2:16}
\mathcal{P}_n([v]_n) \in V_{R+\ep}, \quad \text{for all}~[v]_n\in V_{R}^n \quad \text{and}~ n>n_\ep.
\end{equation}
\end{lemma}

Proof. Let $0<\ep<<R$, $v\in V_{R-\ep}$ and $\mathcal{Q}(v)=[v]_n=([s]_n,[g]_n)$. By applying Cauchy-Bunyakovski-Schwarz (CBS) inequality and Fubini's theorem we have
\begin{gather}
\sum\limits_{k=1}^{n-1}\tau s_{\overline{t}t,k}^2=\sum\limits_{k=1}^{n-1}\frac{1}{\tau^3}\Big [ \int\limits_{t_k}^{t_{k+1}}(s'(t)-s'(t-\tau))dt \Big ]^2 \le \frac{1}{\tau^2} \int_{\tau}^T|s'(t)-s'(t-\tau)|^2 dt \nonumber\\ \le \frac{1}{\tau} \int\limits_{\tau}^{T}dt \int\limits_{t-\tau}^t |s''(\xi)|^2d\xi
\le \int\limits_0^T|s''(t)|^2 dt, \ \sum\limits_{k=1}^{n}\tau s_{\overline{t},k}^2 \le \int\limits_0^T|s'(t)|^2 dt,\label{Eq:W:2:17}
\end{gather}
\begin{gather}
\Big | \sum\limits_{k=0}^{n-1}\tau s_k^2 - \int_0^T s^2(t) dt \Big |= \Big | \sum\limits_{k=0}^{n-1}\int\limits_{t_k}^{t_{k+1}}\int\limits_t^{t_k}(s^2(\xi))'d\xi dt \Big | \le \nonumber\\
\sum\limits_{k=0}^{n-1}\int\limits_{t_k}^{t_{k+1}}\int\limits_{t_k}^t  [ s^2(\xi)+(s'(\xi))^2 \ ]d\xi dt
\le \tau \int\limits_0^T [ s^2(t)+(s'(t))^2 \ ]dt \le R^2 \tau, \label{Eq:W:2:18}
\end{gather}
\begin{gather}
\sum\limits_{k=1}^{n}\tau g_{\overline{t},k}^2 \le \int\limits_0^T|s'(t)|^2 dt, \ 
\Big | \sum\limits_{k=0}^{n-1}\tau g_k^2 - \int_0^T g^2(t) dt \Big | \le R^2 \tau.\label{Eq:W:2:19}
\end{gather}
From (\ref{Eq:W:2:17})-(\ref{Eq:W:2:19}) it follows that
\begin{equation}
\max \Big ( \Vert [s]_n \Vert^2_{w_2^2}, \Vert [g]_n \Vert^2_{w_2^1} \Big ) \le \max \Big ( \Vert s \Vert^2_{W_2^2[0,T]}, \Vert g \Vert^2_{W_2^1[0,T]} \Big )+R^2\tau \le (R-\ep)^2+R^2\tau\le R^2,\label{Eq:W:2:20}
\end{equation}
if $n>n_\ep=\Big [ \frac{RT}{\ep} \Big ] + 1$. Hence, (\ref{Eq:W:2:15}) is proved.

Let us know choose $[v]_n\in V_R^n$. We simplify the notation and assume $v=(s,g)=\mathcal{P}_n([v]_n)$. Through direct calculations we derive
\begin{equation}\label{Eq:W:2:21}
\Vert s \Vert_{W_2^2[0,T]}^2 \le \sum\limits_{k=0}^{n-1}\tau s_k^2 +\sum\limits_{k=1}^{n-1}\tau s_{\overline{t},k}^2 + \sum\limits_{k=1}^{n-1}\tau s_{\overline{t}t,k}^2 + \frac{1}{3}\tau s_{\overline{t},1}^2 + \frac{1}{\tau}s_{\overline{t},1}^2 + C \tau,
\end{equation}
where $C$ is independent of $\tau$. By using CBS inequality we have
\begin{gather}
\tau s_{\overline{t},1}^2 \le \int\limits_0^\tau |s'(t)|^2dt, \ \  \frac{1}{\tau}s_{\overline{t},1}^2
=\frac{1}{\tau^3}\Big | \int\limits_0^\tau \int\limits_0^t s''(\xi)d\xi dt \Big |^2 \le  \nonumber\\ \frac{1}{2\tau}\int\limits_0^\tau \int\limits_0^t |s''(\xi)|^2 d\xi dt \le \frac{1}{2} \int\limits_0^\tau |s''(t)|^2 dt.\label{Eq:W:2:22}
\end{gather}
Since $[v]_n \in V_R^n$, from (\ref{Eq:W:2:21}),(\ref{Eq:W:2:22}) it follows that
\begin{equation}\label{Eq:W:2:23}
\Vert s \Vert_{W_2^2[0,T]}^2 \le C_1,
\end{equation}
where $C_1$ is independent of $\tau$. This implies that
\begin{equation}\label{Eq:W:2:24}
\lim\limits_{\tau \to 0}\Vert s \Vert_{W_2^2[0,\tau]}=0.
\end{equation}
In a similar way we calculate
\begin{equation}\label{Eq:W:2:25}
\Vert g \Vert_{W_2^1[0,T]}^2 \le \sum\limits_{k=0}^{n-1}\tau g_k^2 + \sum\limits_{k=1}^{n}\tau g_{\overline{t},k}^2 + C \tau.
\end{equation}
Hence, from (\ref{Eq:W:2:21}),(\ref{Eq:W:2:22}) and (\ref{Eq:W:2:25}) it follows that
\begin{gather}
\max \Big ( \Vert s \Vert^2_{W_2^2[0,T]}, \Vert g \Vert^2_{W_2^1[0,T]} \Big ) \le \max \Big ( \Vert [s]_n \Vert^2_{w_2^2}, \Vert [g]_n \Vert^2_{w_2^1} \Big )+C\tau 
+\frac{1}{2}\Vert s' \Vert_{W_2^1[0,\tau]}^2 \nonumber\\ \le R^2+C\tau 
+\frac{1}{2}\Vert s' \Vert_{W_2^1[0,\tau]}^2,\label{Eq:W:2:26}
\end{gather}
From (\ref{Eq:W:2:24}) it follows that given $\ep>0$ we can choose $n_{\ep}$ such that for any $n>n_{\ep}$
\begin{equation}\label{Eq:W:2:27}
R^2+C\tau 
+\frac{1}{2}\Vert s' \Vert_{W_2^1[0,\tau]}^2 \le (R+\ep)^2.
\end{equation}
From (\ref{Eq:W:2:26}) and (\ref{Eq:W:2:27}), (\ref{Eq:W:2:16}) follows. Lemma is proved.
\begin{corollary}\label{lipshitz}
Let either $[v]_n \in V_R^n$ or $[v]_n = {\mathcal Q}_n(v)$ for $v \in V_R$. Then
\begin{equation}\label{Eq:W:2:28}
|s_k-s_{k-1}|\le C \tau, \ \ k=1,2,\cdots,n
\end{equation}
where $C$ is independent of $n$. 
\end{corollary}
Indeed, if $v \in V_R$, then $s^{\prime}\in W_2^1[0,T]$ and by Morrey inequality
\begin{equation}\label{Eq:W:2:29}
\Vert s^{\prime} \Vert_{C[0,T]} \le C_1 \Vert s^{\prime} \Vert_{W_2^1[0,T]} \le C_1 R
\end{equation}
and hence for the first component $[s]_n$ of $[v]_n = {\mathcal Q}_n(v)$ we have (\ref{Eq:W:2:28}). Also, if $[v]_n \in V_R^n$, then the sequence $v^n={\mathcal P}_n([v]_n)$ belongs to $V_{R+1}$ by Lemma~\ref{mappings} and 
the component $s^n$ of $v^n$ satisfies (\ref{Eq:W:2:29}). Since,
\[ s^n(0)=s_0, \ s^n(t_k)=\frac{s_k+s_{k-1}}{2}, k=1,\cdots,n \]
from (\ref{Eq:W:2:29}), (\ref{Eq:W:2:28}) easily follows.

\section{Proofs of the Main Results}\label{proofs}
\subsection{First Energy Estimate and its Consequences}\label{firstenergyestimate}
Throughout this section we assume that
\[ \phi \in L_2[0,l], \ \gamma, \chi \in W_2^{1,0}(D), \ a \in L_{\infty}(D),  \]
$a$ satisfies (\ref{Eq:W:1:6}) and $b,c,f$ satisfy the conditions imposed in Section~\ref{E:1:4}.
The main goal of this section to prove the following energy estimation for the discrete state vector.
\begin{theorem}\label{firstenergy}
For all sufficiently small $\tau$ discrete state vector $[u([v]_n)]_n$ satisfies the following stability estimations:
\begin{gather}
\max\limits_{0\le k \le n} \int_{0}^{l} u^2(x;k) \,dx+\tau \sum_{k=1}^{n}\int_0^l \Big | \frac{du(x;k)}{dx} \Big |^2 dx  \le \nonumber\\
C \Big ( \Vert \phi \Vert_{L_2(0,s_0)}^2 + \Vert g \Vert_{L_2(0,T)}^2 + \Vert f \Vert_{L_2(D)}^2 + 
\Vert \gamma(s^n(t),t)(s^n)'(t) \Vert_{L_2(0,T)}^2 \nonumber\\ + \Vert \chi(s^n(t),t) \Vert_{L_2(0,T)}^2 +  \sum_{k=1}^{n-1}{\bf 1}_{+}(s_{k+1}-s_{k})  \int_{s_{k}}^{s_{k+1}} u^2(x;k) dx  \Big ),\label{firstenergyestimate1} 
\end{gather}
\begin{gather}
\max\limits_{0\le k \le n} \int_{0}^{l} u^2(x;k) \,dx+\tau \sum_{k=0}^{n}\int_0^l \Big | \frac{du(x;k)}{dx} \Big |^2 dx + \tau^2 \sum_{k=1}^{n}\int_0^l u_{\overline{t}}^2(x;k)dx \le \nonumber\\
C \Big ( \Vert \phi \Vert_{W_2^1(0,s_0)}^2 + \Vert g \Vert_{L_2(0,T)}^2 + \Vert f \Vert_{L_2(D)}^2 + 
\Vert \gamma(s^n(t),t)(s^n)'(t) \Vert_{L_2(0,T)}^2 \nonumber\\ + \Vert \chi(s^n(t),t) \Vert_{L_2(0,T)}^2 +  \sum_{k=1}^{n-1}{\bf 1}_{+}(s_{k+1}-s_{k})  \int_{s_{k}}^{s_{k+1}} u^2(x;k) dx  \Big ),\label{firstenergyestimate2} 
\end{gather}
where $C$ is independent of $\tau$ and ${\bf 1}_{+}$ be an indicator function of the positive semiaxis.
\end{theorem}
We split the proof into two Lemmas.
\begin{lemma}\label{firstenergy-step1}
For all sufficiently small $\tau$, discrete state vector $[u([v]_n)]_n$ satisfies the following estimation:
\begin{gather}
\max\limits_{1\le k \le n} \int_{0}^{s_k} u^2(x;k) \,dx+\tau \sum_{k=1}^{n}\int_0^{s_k} \Big | \frac{du(x;k)}{dx} \Big |^2 dx + \tau^2 \sum_{k=1}^{n}\int_0^{s_k} u_{\overline{t}}^2(x;k)dx \le \nonumber\\
C \Big ( \Vert \phi \Vert_{L_2(0,s_0)}^2 + \Vert g \Vert_{L_2(0,T)}^2 + \Vert f \Vert_{L_2(D)}^2 + 
\Vert \gamma(s^n(t),t)(s^n)'(t) \Vert_{L_2(0,T)}^2 \nonumber\\ + \Vert \chi(s^n(t),t) \Vert_{L_2(0,T)}^2 +  \sum_{k=1}^{n-1}{\bf 1}_{+}(s_{k+1}-s_{k}) \int_{s_{k}}^{s_{k+1}} u^2(x;k) dx \Big ),\label{firstenergyestimate3} 
\end{gather}
where $C$ is independent of $\tau$.
\end{lemma}

Proof. By choosing $\eta(x)=2\tau u(x;k)$ in (\ref{Eq:W:1:13}) and by using the equality
\[ 2\tau u_{\overline{t}}(x;k)u(x;k) = u^2(x;k)-u^2(x;k-1)+\tau^2 u_{\overline{t}}^2(x;k)  \]
we have
\begin{gather}
\int_{0}^{s_{k}}u^2(x;k)dx-\int_{0}^{s_{k}}u^2(x;k-1)dx+\tau^2\int_{0}^{s_{k}}u_{\overline{t}}^2(x;k)dx+ 2\tau\int_{0}^{s_{k}}a_{k}(x)\Big |\frac{d u(x;k)}{d x}\Big |^2 dx = \nonumber\\
2\tau\int_{0}^{s_{k}} b_{k}(x)\frac{d u(x;k)}{d x}u(x;k)+ c_{k}(x)u^2(x;k)
 - f_{k}(x) u(x;k) \Big]\, dx - \nonumber\\
 2 \tau \left[ (\gamma_{s^n} (s^n)')^k-\chi^{k}_{s^n} \right] u(s_{k};k)
 -2\tau g^{k}u(0;k).\label{Eq:W:3:4}
\end{gather}
Using (\ref{Eq:W:1:6}), Cauchy inequalities with appropriately chosen  $\ep >0$, and Morrey inequality \eqref{Morrey} from (\ref{Eq:W:3:4}) we derive that
\begin{gather}
\int_{0}^{s_{k}}u^2(x;k)dx-\int_{0}^{s_{k}}u^2(x;k-1)dx+ a_0\tau\int_{0}^{s_{k}}\Big |\frac{d u(x;k)}{d x}\Big |^2 dx +   \tau^2\int_{0}^{s_{k}}u_{\overline{t}}^2(x;k)dx \le \nonumber\\
 C_1 \tau \left[ |(\gamma_{s^n} (s^n)')^k|^2+|\chi^{k}_{s^n}|^2+|g^k|^2+ \int_{0}^{s_{k}}f_k^2(x)dx+\int_{0}^{s_{k}}u^2(x;k)dx\right],\label{Eq:W:3:5}
\end{gather}
where $C_1$ is independent of $\tau$. Assuming that$\tau < C_1$, from (\ref{Eq:W:3:5}) it follows that
\begin{gather}
(1-C_1\tau)\int_{0}^{s_{k}}u^2(x;k)dx \le 
\int_{0}^{s_{k-1}}u^2(x;k-1)dx+ {\bf 1}_{+}(s_k-s_{k-1})\int_{s_{k-1}}^{s_k}u^2(x;k-1)dx + \nonumber\\ 
 C_1 \tau \left[ |(\gamma_{s^n} (s^n)')^k|^2+|\chi^{k}_{s^n}|^2+|g^k|^2+ \int_{0}^{s_{k}}f_k^2(x)dx\right],\label{Eq:W:3:6}
\end{gather}
By induction we have
\begin{gather}
\int_{0}^{s_{k}}u^2(x;k)dx \le (1-C_1\tau)^{-k}\int_0^{s_0}\phi^2(x)dx +  \sum_{j=1}^{k} (1-C_1\tau)^{-k+j-1} \Big \{ C_1 \tau \times \nonumber\\ 
 \Big [ |(\gamma_{s^n} (s^n)')^j|^2+|\chi^{j}_{s^n}|^2+|g^j|^2+ \int_{0}^{s_{j}}f_j^2(x)dx \Big ]+{\bf 1}_{+}(s_j-s_{j-1})\int_{s_{j-1}}^{s_j}u^2(x;j-1)dx \Big \}.\label{Eq:W:3:7}
\end{gather}
For arbitrary $1\le j \le k \le n$ we have
\begin{equation}\label{Eq:W:3.8}
(1-C_1\tau)^{-k+j-1}\le(1-C_1\tau)^{-k}\le(1-C_1\tau)^{-n}=\Big ( 1-\frac{C_1 T}{n}\Big )^{-n} \to e^{C_1T},
\end{equation}
as $\tau \to 0$. Accordingly for sufficiently small $\tau$ we have
\begin{equation}\label{Eq:W:3:9}
(1-C_1\tau)^{-k+j-1}\le 2e^{C_1T} \  \quad \text{for} \ 1\le j \le k \le n,
\end{equation}
By applying CBS inequality from (\ref{Eq:W:3:7})-(\ref{Eq:W:3:9}) it follows that
\begin{gather}
\max\limits_{1\le k \le n} \int_{0}^{s_k} u^2(x;k) \,dx \le 
C_2 \Big ( \Vert \phi \Vert_{L_2(0,s_0)}^2 + \Vert g \Vert_{L_2(0,T)}^2  + 
\Vert \gamma(s^n(t),t)(s^n)'(t) \Vert_{L_2(0,T)}^2 + \nonumber\\ \Vert \chi(s^n(t),t) \Vert_{L_2(0,T)}^2 + \Vert f \Vert_{L_2(D)}^2 +  \sum_{k=1}^{n-1} {\bf 1}_{+}(s_{k+1}-s_{k})\int_{s_{k}}^{s_{k+1}} u^2(x;k) dx  \Big ).\label{Eq:W:3:10} 
\end{gather}
where $C_2$ is independent of $\tau$. Having (\ref{Eq:W:3:10}), we perform summation of (\ref{Eq:W:3:5})
with respect to $k$ from $1$ to $n$ and derive 
\begin{gather}
\int_{0}^{s_n} u^2(x;n) \,dx+a_0\tau \sum_{k=1}^{n}\int_0^{s_k} \Big | \frac{du(x;k)}{dx} \Big |^2 dx + \tau^2 \sum_{k=1}^{n}\int_0^{s_k} u_{\overline{t}}^2(x;k)dx \le \nonumber\\
\Vert \phi \Vert_{L_2(0,s_0)}^2 + C_1 \Big ( \Vert g \Vert_{L_2(0,T)}^2 + \Vert f \Vert_{L_2(D)}^2 + 
\Vert \gamma(s^n(t),t)(s^n)'(t) \Vert_{L_2(0,T)}^2 \nonumber\\ + \Vert \chi(s^n(t),t) \Vert_{L_2(0,T)}^2 +  \tau \sum_{k=1}^{n}\int_{0}^{s_{k}}u^2(x;k)dx \Big ) +\sum_{k=1}^{n-1}{\bf 1}_{+}(s_{k+1}-s_k) \int_{s_{k}}^{s_{k+1}} u^2(x;k) dx,\label{Eq:W:3:11} 
\end{gather}
From (\ref{Eq:W:3:10}) and (\ref{Eq:W:3:11}), (\ref{firstenergyestimate3}) follows. Lemma is proved.

In the next lemma we prove a nice property of the extension introduced in the Definition 1.3, which allows
to extend the estimation (\ref{firstenergyestimate3}) to (\ref{firstenergyestimate1}) and (\ref{firstenergyestimate2}).
\begin{lemma}\label{extensionlemma}
Given discrete control vector $[v]_n \in V_R^n$, a discrete state vector $[u([v]_n)]_n$ satisfies the inequalty
\begin{gather}
\max\limits_{1\le k \le n} \int_{0}^{l} u^2(x;k) \,dx+\tau \sum_{k=0}^{n}\int_0^l \Big | \frac{du(x;k)}{dx} \Big |^2 dx + \tau \sum_{k=1}^{n}\int_0^l u_{\overline{t}}^2(x;k)dx \le \nonumber\\
C \Big ( \max\limits_{1\le k \le n} \int_{0}^{s_k} u^2(x;k) \,dx+\tau \sum_{k=0}^{n}\int_0^{s_k} \Big | \frac{du(x;k)}{dx} \Big |^2 dx + \tau \sum_{k=1}^{n}\int_0^{s_k} u_{\overline{t}}^2(x;k)dx \Big ),\label{Eq:W:3:12}
\end{gather}
where $C$ is independent of $\tau$. 
\end{lemma}

Proof. By induction it follows that the first two terms on the left hand side are estimated by the first two terms on the right hand side with the constant $C=2^N$, where  $N$ is defined in
(\ref{Eq:W:1:14}). 


Define a family of functions $\{\tilde{u}(y;k),k=0,...,n\}$ as
\[ \tilde{u}(y;0)=\phi(ys_0), \  \tilde{u}(y;k)=u(ys_k;k), \ \ 0\le y \le 1, k=1,...,n. \]
As before, assume they are all continued by induction to semiaxis $\{y\ge 0\}$  as
\[ \tilde{u}(y;k)=\tilde{u}(2^n-y;k), \quad \text{for}~ 2^{n-1} \le y \le 2^n. \] 
We have
\begin{gather}
\sum\limits_{k=1}^n \tau \int\limits_0^l u_{\overline{t}}^2(x;k)dx \le \sum\limits_{k=1}^n \tau \int\limits_0^{2^Ns_k}u_{\overline{t}}^2(x;k)dx = \nonumber\\ \sum\limits_{k=1}^n \tau \int\limits_0^{2^Ns_k}\Big [  \frac{\tilde{u}(x/s_k; k)-\tilde{u}(x/s_{k-1}; k-1)}{\tau}\Big ]^2 dx = 
\nonumber\\ \sum\limits_{k=1}^n \tau s_k \int\limits_0^{2^N}\Big [  \frac{\tilde{u}(y; k)-\tilde{u}(y s_k/s_{k-1}; k-1)}{\tau}\Big ]^2 dy \le I_1+I_2 \label{Eq:W:3:13}
\end{gather}
where
\begin{gather}
I_1=2 \sum\limits_{k=1}^n \tau s_k \int\limits_0^{2^N}\Big [  \frac{\tilde{u}(y; k)-\tilde{u}(y; k-1)}{\tau}\Big ]^2 dy=\cdots=  2^{N+1}\sum\limits_{k=1}^n \tau s_k \int\limits_0^1 \tilde{u}_{\overline{t}}^2(y;k)dy= \nonumber\\ 2^{N+1}\sum\limits_{k=1}^n \tau \int\limits_0^{s_k}\Big [ \frac{ u(x; k)-u(x s_{k-1}/s_k; k-1)}{\tau}\Big ]^2 dx  
\le \nonumber\\ 2^{N+2}\sum\limits_{k=1}^n \tau \int\limits_0^{s_k} u_{\overline{t}}^2(x;k)dx + 2^{N+2}\sum\limits_{k=1}^n \tau \int\limits_0^{s_k}\Big [ \frac{  u(x; k-1)-u(x s_{k-1}/s_k; k-1)}{\tau}\Big ]^2 dx   
\label{Eq:W:3:14}
\end{gather}
\begin{equation}\label{Eq:W:3:15}
I_2=2 \sum\limits_{k=1}^n \tau s_k \int\limits_0^{2^N}\Big [  \frac{\tilde{u}(y; k-1)-\tilde{u}(ys_k/s_{k-1}; k-1)}{\tau}\Big ]^2 dy.
\end{equation}
By using CBS inequality, Fubini's theorem and Corollary~\ref{lipshitz} we have
\begin{gather}
\sum\limits_{k=1}^n \tau \int\limits_0^{s_k}\Big [ \frac{  u(x; k-1)-u(x s_{k-1}/s_k; k-1)}{\tau}\Big ]^2 dx = \nonumber\\ \sum\limits_{k=1}^n \frac{1}{\tau}\int\limits_0^{s_k}\Big | \int\limits_{x \frac{s_{k-1}}{s_k}}^x \frac{du(\xi;k-1)}{d\xi} d\xi \Big |^2dx  \le \frac{C_1^2 l}{\delta}\sum\limits_{k=0}^{n-1} \tau \int\limits_0^l \Big | \frac{du(x;k)}{dx} \Big |^2dx,\label{Eq:W:3:16}
\end{gather}
\begin{gather}
I_2\le \frac{2^{2N+1}C_1^2N^2}{\delta}\sum\limits_{k=1}^n \tau \int\limits_0^{N2^N} \Big |\frac{d\tilde{u}(x;k-1)}{dx} \Big |^2 dx = \frac{2^{3N+1}C_1^2N^3}{\delta}\sum\limits_{k=0}^{n-1} \tau \int\limits_0^{1} \Big |\frac{d\tilde{u}(x;k)}{dx} \Big |^2 dx \nonumber\\
\le  \frac{2^{3N+1}C_1^2N^3l}{\delta}\sum\limits_{k=0}^{n-1} \tau \int\limits_0^{s_k} \Big |\frac{du(x;k)}{dx} \Big |^2 dx \label{Eq:W:3:17} 
\end{gather}
Hence, from (\ref{Eq:W:3:13})-(\ref{Eq:W:3:17}) it follows that 
\begin{equation}\label{Eq:W:3:18}
\sum\limits_{k=1}^n \tau \int\limits_0^l u_{\overline{t}}^2(x;k)dx \le C \Big (\sum\limits_{k=0}^{n-1} \tau \int\limits_0^{s_k} \Big |\frac{du(x;k)}{dx} \Big |^2 dx + \sum\limits_{k=1}^n \tau \int\limits_0^{s_k} u_{\overline{t}}^2(x;k)dx \Big )
\end{equation}
where $C$ is independent of $\tau$. From (\ref{Eq:W:3:18}),(\ref{Eq:W:3:12}) follows. Lemma is proved.

It can be easily seen that Theorem~\ref{firstenergy} follows from Lemma~\ref{firstenergy-step1} and
Lemma~\ref{extensionlemma}.

Let $[v]_n \in V_R^n, n=1,2,...$ be a sequence of discrete controls. From  Lemma~\ref{mappings}
it follows that the sequence $\{\mathcal{P}_n([v]_n)\}$ is weakly precompact in $W_2^2[0,T]\times W_2^1[0,T]$. Assume that the whole sequence converges to $v=(s,g)$
weakly in $W_2^2[0,T]\times W_2^1[0,T]$. This implies the strong convegence in $W_2^1[0,T] \times L_2[0,T]$. Conversely, given control $v=(s,g)\in V_R^n$ we can choose sequence of discrete
controls $[v]_n = \mathcal{Q}_n(v)$. Appplying Lemma~\ref{mappings} twice one can easily establish
that the sequence $\{\mathcal{P}_n([v]_n\}$ converges to $v=(s,g)$
weakly in $W_2^2[0,T]\times W_2^1[0,T]$, and strongly in $W_2^1[0,T]\times L_2[0,T]$. In the next theorem we prove the continuous dependence of the family of interpolarions $\{u^\tau\}$ on this convergence.
\begin{theorem}\label{continuity1}
Let $[v]_n \in V_R^n, n=1,2,...$ be a sequence of discrete controls and the sequence $\{\mathcal{P}_n([v]_n\}$ converges strongly in $W_2^1[0,T] \times L_2[0,T]$ to $v=(s,g)$. Then the sequence $\{u^\tau\}$ converges as $\tau \to 0$ weakly in $W_2^{1,0}(\Omega)$ to weak solution $u \in V_2^{1,0}(\Omega)$ of the problem (\ref{Eq:W:1:1})-(\ref{Eq:W:1:4}), i.e. to the solution of the integral identity (\ref{Eq:W:1:10}). Moreover, $u$ satisfies the energy estimate
\begin{equation}\label{V210estimate}
\Vert u \Vert_{V_2^{1,0}(D)}^2 \le C \Big ( \Vert \phi \Vert_{L_2(0,s_0)}^2 + \Vert g \Vert_{L_2(0,T)}^2 + \Vert f \Vert_{L_2(D)}^2 + 
\Vert \gamma \Vert_{W_2^{1,0}(D)}^2 + \Vert \chi \Vert_{W_2^{1,0}(D)}^2 \Big )
\end{equation}
\end{theorem}

Proof. In addition to quadratic interpolation of $[s]_n$ from (\ref{Eq:W:1:11}), consider two linear interpolations:
\begin{equation*}
\tilde{s}^n(t)=s_{k-1}+\frac{s_k-s_{k-1}}{\tau}(t-t_{k-1}), \ \ t_{k-1} \le t \le t_k, k=\overline{1,n}; \ \tilde{s}^n(t) \equiv s_n, \ \ t\ge T;
\end{equation*}
\begin{equation*}
\tilde{s}_1^n(t)=\tilde{s}^n(t+\tau), \ \ 0\le t \le T.
\end{equation*}
It can be easily proved that both sequences $\tilde{s}^n$ and $\tilde{s}_1^n$ are equivalent 
to the sequence $s^n$ in $W_2^1[0,T]$ and converge to $s$ strongly in $W_2^1[0,T]$. In particular,
\begin{equation}\label{Eq:W:3:19}
\sup\limits_{n}\Vert \tilde{s}_1^n \Vert_{W_2^1[0,T]} < C_*
\end{equation}
where $C_*$ is independent of $n$. We estimate the last term on the right-hand side of (\ref{firstenergyestimate1})
as follows:
\begin{gather}
\sum_{k=1}^{n-1}{\bf 1}_{+}(s_{k+1}-s_{k})  \int_{s_{k}}^{s_{k+1}} u^2(x;k) dx =
\sum_{k=1}^{n-1}{\bf 1}_{+}(s_{k+1}-s_{k})  \int_{t_{k}}^{t_{k+1}}(\tilde{s}^n)^{'}(t) u^2(\tilde{s}^n(t);k) dt = \nonumber\\ \sum_{k=1}^{n-1}{\bf 1}_{+}(s_{k+1}-s_{k})  \int_{t_{k}}^{t_{k+1}}(\tilde{s}^n)^{'}(t) \Big ( u^{\tau}(\tilde{s}^n(t),t-\tau) \Big )^2 dt = \nonumber\\ \sum_{k=1}^{n-1}{\bf 1}_{+}(s_{k+1}-s_{k})  \int_{t_{k-1}}^{t_{k}}(\tilde{s}_1^n)^{'}(t) \Big ( u^{\tau}(\tilde{s}_1^n(t),t) \Big )^2 dt. \label{Eq:W:3:20}
\end{gather}
By applying CBS inequality we have
\begin{equation}\label{Eq:W:3:21}
\Big | \sum_{k=1}^{n-1}{\bf 1}_{+}(s_{k+1}-s_{k})  \int_{s_{k}}^{s_{k+1}} u^2(x;k) dx \Big | \le
\Vert (\tilde{s}_1^n)^{\prime}\Vert_{L_2[0,T]} \Vert u^{\tau}(\tilde{s}_1^n(t),t)\Vert_{L_4[0,T]}^2.
\end{equation}
From the results on traces of the elements of space $V_2(D)$ (\cite{LSU,BIN,Nikolski}) it follows that for arbitrary $u\in V_2(D)$ the following inequality is valid
\begin{equation}\label{Eq:W:3:22}
\Vert u(\tilde{s}_1^n(t),t)\Vert_{L_4[0,T]}\le \tilde{C}\Vert u \Vert_{V_2(D)},
\end{equation}
with the constant $\tilde{C}$ being independent of $u$ as well as $n$. From (\ref{Eq:W:3:19}),(\ref{Eq:W:3:21}) and (\ref{Eq:W:3:22}) it follows that
\begin{equation}\label{Eq:W:3:23}
\Big | \sum_{k=1}^{n-1}{\bf 1}_{+}(s_{k+1}-s_{k})  \int_{s_{k}}^{s_{k+1}} u^2(x;k) dx \Big | \le
C_* \tilde{C} \Vert u^{\tau} \Vert_{V_2(D)}^2.
\end{equation}
If the constant $C_*$ from (\ref{Eq:W:3:19}) satisfies the condition
\begin{equation}\label{Eq:W:3:24}
C_*<(C\tilde{C})^{-1}
\end{equation}
then from (\ref{firstenergyestimate1}) and (\ref{Eq:W:3:23}) it follows that
\begin{gather}
\Vert u^{\tau} \Vert_{V_2^{1,0}(D)}^2 \le C \Big ( \Vert \phi \Vert_{L_2(0,s_0)}^2 + \Vert g \Vert_{L_2(0,T)}^2 + \Vert f \Vert_{L_2(D)}^2 + \nonumber\\
\Vert \gamma(s^n(t),t)(s^n)'(t) \Vert_{L_2(0,T)}^2 + \Vert \chi(s^n(t),t) \Vert_{L_2(0,T)}^2 \Big ),\label{Eq:W:3:25} 
\end{gather}
where $C$ is another constant independent of $n$. By applying the results on the traces of elements of $W_2^{1,0}(D)$ (\cite{BIN,Nikolski}) on smooth curve $x=s^n(t)$, Morrey inequality for $(s^n)^{\prime}$ and (\ref{Eq:W:2:16}) we have
\begin{gather}
\Vert \gamma(s^n(t),t)(s^n)'(t) \Vert_{L_2(0,T)} \le \Vert (s^n)^{\prime} \Vert_{C[0,T]} \Vert \gamma(s^n(t),t) \Vert_{L_2[0,T]} \le C_3 \Vert \gamma \Vert_{W_2^{1,0}(D)}\nonumber\\
\Vert \chi(s^n(t),t)\Vert_{L_2[0,T]} \le C_3 \Vert \chi \Vert_{W_2^{1,0}(D)}, \label{Eq:W:3:26}
\end{gather}
where $C_3$ is independent of $\gamma, \chi$ and $n$. Hence, from (\ref{Eq:W:3:25}) and (\ref{Eq:W:3:26})
it follows the estimation
\begin{equation}\label{Eq:W:3:27}
\Vert u^{\tau} \Vert_{V_2^{1,0}(D)}^2 \le C \Big ( \Vert \phi \Vert_{L_2(0,s_0)}^2 + \Vert g \Vert_{L_2(0,T)}^2 + \Vert f \Vert_{L_2(D)}^2 + 
\Vert \gamma \Vert_{W_2^{1,0}(D)}^2 + \Vert \chi \Vert_{W_2^{1,0}(D)}^2 \Big ),
\end{equation}
with $C$ being independent of $n$.


If (\ref{Eq:W:3:24}) is not satisfied, then we can partition $[0,T]$ into finitely many segments $[t_{n_{j-1}},t_{n_{j}}]$, $j=\overline{1,q}$ with $t_{n_{0}}=0$, $t_{n_{q}}=T$ in such a way that by replacing $[0,T]$ with any of the subsegments $[t_{n_{j-1}},t_{n_{j}}]$ (\ref{Eq:W:3:19}) will be satisfied with $C_*$ small enough to obey (\ref{Eq:W:3:24}). Hence, we divide $D$ into finitely many subsets
\[ D^{j}=D \cap \{ t_{n_{j-1}}<t\leq t_{n_{j}}\}\]
such that every norm $\Vert u^{\tau}\Vert_{V_{2}(D^{j})}^2$ is uniformly bounded through the right-hand side of (\ref{Eq:W:3:27}). Summation with $j=1,\ldots,q$ implies (\ref{Eq:W:3:27}).

From (\ref{Eq:W:3:27}) it follows that the sequence $\{ u^{\tau}\}$ is weakly precompact in $W_{2}^{1,0}(D)$.  Let $u \in W_{2}^{1,0}(D)$ be a weak limit point of $u^{\tau}$ in $W_{2}^{1,0}(D)$, and assume that whole sequence $\{ u^{\tau}\}$ converges to $u$ weakly in $W_{2}^{1,0}(D)$.  Let us prove that in fact $u$ satisfies the integral identity (\ref{Eq:W:1:10}) for arbitrary test function $\Phi \in W_{2}^{1,1}(\Omega)$ such that $\left. \Phi\right|_{t=T}=0$. Due to density of $C^1(\overline{\Omega})$ in $W_{2}^{1,1}(\Omega)$ it is enough to assume $\Phi \in C^1(\overline{\Omega})$. Without loss of generality we can also assume that
$\Phi \in C^1(\overline{D}_{T+\tau}), \ \Phi\equiv 0, \quad\text{for}~ T\le t \le T+\tau$, where
\[ D_{T+\tau}=\{ (x,t):~0<x<l,~0<t\leq T+\tau\} \]
Otherwise, we can continue $\Phi$ to $D_{T+\tau}$ with the described properties.

Let 
\[\Phi(x;k)=\Phi(x,k\tau),\quad \Phi_{t}(x;k)=\frac{\Phi(x;k+1)-\Phi(x;k)}{\tau}\]
As before, we construct piecewise constant interpolations $\Phi^{\tau}$, $\Phi_{t}^{\tau}$.  Obviously, the sequences $\{ \Phi^{\tau}\}$, $\{ \frac{\partial \Phi^{\tau}}{\partial x}\}$ and $\{ \Phi_{t}^{\tau}\}$ converge as $\tau \to 0$ uniformly in $\overline{D}$ to $\Phi$, $\frac{\partial \Phi}{\partial x}$ and $\frac{\partial \Phi}{\partial t}$ respectively.

By choosing in (\ref{Eq:W:1:13})
 $\eta(x)=\tau\Phi(x;k)$, after summation with respect to $k=\overline{1,n}$ and transformation of the time difference term as follows
\begin{gather}
\tau\sum_{k=1}^{n}\int_{0}^{s_{k}}u_{\bar{t}}(x;k)\Phi(x;k)\, dx=
-\tau \sum_{k=1}^{n-1}\int_{0}^{s_{k+1}}u(x;k)\Phi_{t}(x;k)\, dx- \int_{0}^{s_{1}}\phi(x)\Phi(x;1)\, dx \nonumber\\ -\sum_{k=1}^{n-1}\int_{s_{k}}^{s_{k+1}}u(x;k)\Phi(x;k)\, dx=
-\sum_{k=1}^{n}\int_{t_{k-1}}^{t_{k}}\int_{0}^{s_{k+1}}u^{\tau}\Phi_{t}^{\tau}\, dx\, dt
-\int_{0}^{s_{1}}\phi(x)\Phi^{\tau}(x,\tau)\, dx \nonumber\\
- \sum_{k=1}^{n-1}\int_{t_{k}}^{t_{k+1}}(\tilde{s}^{n})^{\prime}(t)u^\tau(\tilde{s}^{n}(t),t-\tau)\Phi^\tau(\tilde{s}^{n}(t),t-\tau) \, dt = -\int_{0}^{T}\int_{0}^{s(t)}u^{\tau}\Phi_{t}^{\tau}\, dx \, dt\nonumber\\ -\int_{0}^{s_{1}}\phi(x)\Phi^{\tau}(x,\tau)\, dx-\int_{0}^{T-\tau}(\tilde{s}^{n}_1)^{\prime}(t)u^{\tau}((\tilde{s}^{n}_1)(t),t)\Phi^{\tau}((\tilde{s}^{n}_1)(t),t)\, dt\nonumber\\
-\sum_{k=1}^{n-1}\int_{t_{k-1}}^{t_{k}}\int_{s(t)}^{s_{k+1}}u^{\tau}\Phi_{t}^{\tau}\, dx\, dt\label{Eq:W:3:28}
\end{gather}
we derive that

\begin{gather}
\int_{0}^{T}\int_{0}^{s(t)}\bigg\{a \frac{\partial u^{\tau}}{\partial x}\frac{\partial \Phi^{\tau}}{\partial x}  - b \frac{\partial u^{\tau}}{\partial x}\Phi^{\tau}-c u^{\tau}\Phi^{\tau} + f \Phi^{\tau} -u^{\tau}\Phi_{t}^{\tau}\bigg\}\, dx \,dt-\int_{0}^{s_0}\phi(x)\Phi^{\tau}(x,\tau)\, dx\nonumber\\
-\int_{0}^{T-\tau}(\tilde{s}^{n}_1)^{\prime}(t)u^{\tau}((\tilde{s}^{n}_1)(t),t)\Phi^{\tau}((\tilde{s}^{n}_1)(t),t)\, dt
+\int_{0}^{T}g(t)\Phi^{\tau}(0,t)\, dt\nonumber\\
+\int_{0}^{T}\Big[\gamma(s^n(t),t)(s^n)^{\prime}(t)- \chi(s^n(t),t)) \Big]\Phi^{\tau}(s^n(t),t)\, dt
-R=0\label{Eq:W:3:29}
\intertext{where}
R=\sum_{k=1}^{n}\int_{t_{k-1}}^{t_{k}}\int_{s(t)}^{s_{k}}\bigg\{a \frac{\partial u^{\tau}}{\partial x}\frac{\partial \Phi^{\tau}}{\partial x}  - b \frac{\partial u^{\tau}}{\partial x}\Phi^{\tau}-c u^{\tau}\Phi^{\tau} + f \Phi^{\tau}\bigg\}\, dx \,dt - \sum_{k=1}^{n-1}\int_{t_{k-1}}^{t_{k}}\int_{s(t)}^{s_{k+1}}u^{\tau}\Phi_{t}^{\tau}\, dx\, dt\nonumber\\
+\sum_{k=1}^{n}\int_{t_{k-1}}^{t_{k}}\int_{s^n(t)}^{s_{k}}\Big[ \gamma(s^n(t),t) (s^n)^{\prime}(t) - \chi(s^n(t),t))\Big]\frac{\partial \Phi^{\tau}}{\partial x}\, dx \, dt + \int_{s_0}^{s_1}\phi(x)\Phi^{\tau}(x,\tau)\, dx \nonumber
\end{gather}
Let
\[ \Delta = \bigcup_{k=1}^{n}\left\{ (x,t):t_{k-1}<t<t_{k},~\min(s(t),s_{k}) < x < \max(s(t),s_{k})\right\}\]
$|\Delta|$ denotes the Lebesgue measure of $\Delta$.
Since
\[ |\Delta| \leq \sum_{k=1}^{n}\int_{t_{k-1}}^{t_{k}}\int_{t}^{t_{k}}|s'(\tau)|\, d\tau \, dt
\leq \frac{2\sqrt{T}}{3}\Vert s'\Vert_{L^{2}(0,T)}\tau \to 0\quad \text{as}~\tau \to 0 \]
and all of the integrands are uniformly bounded in $L^{1}(D)$, it follows that the first term in the expression of $R$ converges to zero as $\tau \to 0$. In a similar way one can see that the second and third terms also converge to zero as $\tau \to 0$. The last term in the expression of $R$ converges to zero due to Corollary~\ref{lipshitz} and uniform convergence of $\{\Phi^\tau \}$ in $\overline{D}$. Hence, we have
\begin{equation}
\lim_{\tau \to 0}R=0\label{Eq:W:3:30}
\end{equation}
Due to weak convergence of $u^\tau$ to $u$ in $W_2^{1,0}(D)$ and uniform convergence of the sequences $\{ \Phi^{\tau}\}$, $\{ \frac{\partial \Phi^{\tau}}{\partial x}\}$ and $\{ \Phi_{t}^{\tau}\}$  to $\Phi$, $\frac{\partial \Phi}{\partial x}$ and $\frac{\partial \Phi}{\partial t}$ respectively, passing to limit as $\tau \to 0$, it follows that first, second and fourth integrals on the left-hand side of (\ref{Eq:W:3:29}) converge to similar integrals with $u^\tau$, $\Phi^\tau$, $\Phi^\tau_t$, $\Phi^\tau(x,\tau)$ and $\Phi^\tau(0,t)$ replaced by $u$,$\Phi$, $\frac{\partial \Phi}{\partial t}$, $\Phi(x,0)$ and $\Phi(0,t)$ respectively. Since $s^n$ converges to $s$ strongly in $W_2^1[0,T]$, the traces $\gamma(s^n(t),(t))$, $\chi(s^n(t),t)$ converge strongly in $L_2[0,T]$ to traces $\gamma(s(t),(t))$, $\chi(s(t),t)$ respectively. Since $\Phi^{\tau}(s^n(t),t)$  converge uniformly on $[0,T]$ to
$\Phi(s(t),t)$, passing to the limit as $\tau \to 0$, the last integral on the left-hand side of 
(\ref{Eq:W:3:29}) converge to similar integral with $s^n$ and $\Phi^\tau$ replaced by $s$ and $\Phi$.

It only remains to prove that
\begin{equation}
\lim_{\tau \to 0}\int_{0}^{T-\tau}(\tilde{s}^{n}_1)^{\prime}(t)u^{\tau}(\tilde{s}^{n}_1(t),t)\Phi^{\tau}(\tilde{s}^{n}_1(t),t)\, dt
=\int_{0}^{T}s'(t)u(s(t),t)\Phi(s(t),t)\label{Eq:W:3:31}
\end{equation}
Since $\{\tilde{s}^{n}_1\}$ converges to $s$ strongly in $W_2^1[0,T]$, from (\ref{Eq:W:3:27}) it follows that $\{ u^{\tau}(\tilde{s}^{n}_1(t),t)\}$ is uniformly bounded in $L_{2}[0,T]$ and
\begin{equation}\label{Eq:W:3:32}
\Vert u^{\tau}(\tilde{s}^{n}_1(t),t) - u^{\tau}(s(t),t)\Vert_{L_2[0,T]} \to 0\quad \text{as}~\tau \to 0
\end{equation}
Since $\{u^{\tau}\}$ converges to $u$ weakly in $W_2^{1,0}(D)$, it follows that
\begin{equation}
u^{\tau}(s(t),t) \to u(s(t),t),\quad \text{weakly in}~L_{2}[0,T]\label{Eq:W:3:33}
\end{equation}
Since $\{\Phi^{\tau}(\tilde{s}^{n}_1(t),t)\}$ converges to $\Phi(s(t),t)$ uniformly in $[0,T]$, from
(\ref{Eq:W:3:32}),(\ref{Eq:W:3:33}), (\ref{Eq:W:3:31}) easily follows.

Passing to the limit as $\tau \to 0$, from (\ref{Eq:W:3:29}) it follows that $u$ satisfies integral identity (\ref{Eq:W:1:10}), i.e. it is a weak solution of the problem (\ref{Eq:W:1:1})-(\ref{Eq:W:1:4}). Since this solution is unique (\cite{LSU}) it follows that indeed the whole sequence $\{ u^{\tau}\}$ converges to $u \in V_2^{1,0}(\Omega)$ weakly in $W_2^{1,0}(\Omega)$. From the property of weak convergence and (\ref{Eq:W:3:27}),(\ref{V210estimate}) follows. Lemma is proved.

In particular, Theorem~\ref{continuity1} implies the following well-known existence result (\cite{LSU}):
\begin{corollary}\label{V210solution}
For arbitrary $v=(s,g)\in V_R$ there exists a weak solution $u \in V_2^{1,0}(\Omega)$ of the problem (\ref{Eq:W:1:1})-(\ref{Eq:W:1:4}) which satisfy the energy estimate (\ref{V210estimate})
\end{corollary}
\subsection{Second Energy Estimate and its Consequences}\label{secondenergy}
The main goal of this section to prove the following energy estimation for the discrete state vector.
\begin{theorem}\label{secondenergyestimate}
For all sufficiently small $\tau$ discrete state vector $[u([v]_n)]_n$ satisfies the following stability estimation:
\begin{gather}
\max_{1 \leq k \leq n}\int_{0}^{l}\left| \frac{d u(x;k)}{d x}\right|^{2}\, dx+\tau \sum_{k=1}^{n}\int_{0}^{l}u_{\bar{t}}^{2}(x;k)\, dx  \le \nonumber\\ C \bigg[ \left\Vert \phi\right\Vert_{W_{2}^{1}[0,l]}^{2}+\left\Vert g\right\Vert_{W_{2}^{\frac{1}{4}}[0,T]}^{2}+\left\Vert \gamma(s^n(t),t)(s^n)'(t)\right\Vert_{W_{2}^{\frac{1}{4}}[0,T]}^{2}+\left\Vert \chi(s^n(t),t)\right\Vert_{W_{2}^{\frac{1}{4}}[0,T]}^{2}\nonumber+\\
+\left\Vert f\right\Vert_{L_{2}(D)}^{2} +  \sum_{k=1}^{n-1}{\bf 1}_{+}(s_{k+1}-s_{k})  \int_{s_{k}}^{s_{k+1}} u^2(x;k) dx \bigg] ,\label{secondenergyestimate1} 
\end{gather}
\end{theorem}
We split the proof into two lemmas.
\begin{lemma}
Let given discrete control vector $[v]_n$, along with discrete state vector $[u([v]_n)]_n$, the vector function 
\[ [\tilde{u}([v]_n)]_n=(\tilde{u}(x;0),\tilde{u}(x;1),...,\tilde{u}(x;n))  \]
is defined as
\[
\tilde{u}(x;k)=
\left\{
\begin{array}{l}
u(x;k) \ \ 0\le x \le s_k,\\
u(s_k;k) \ \ s_k\le x \le l, k=\overline{0,n}.
\end{array}\right.
\]
Then for all sufficiently small $\tau$, $[\tilde{u}([v]_n)]_n$ satisfies the following estimation:
\begin{gather}
\max_{1 \leq k \leq n}\int_{0}^{s_{k}}\left| \frac{d \tilde{u}(x;k)}{d x}\right|^{2}\, dx+\tau \sum_{k=1}^{m}\int_{0}^{s_{k}}\tilde{u}_{\bar{t}}^{2}(x;k)\, dx + \tau^{2}\sum_{k=1}^{m}\int_{0}^{s_{k}}\left[ \left( \frac{d \tilde{u}(x;k)}{d x}\right)_{\bar{t}}\right]^{2}\,dx\leq
 \nonumber\\ C \bigg[ \left\Vert \phi\right\Vert_{W_{2}^{1}[0,l]}^{2}+\left\Vert g\right\Vert_{W_{2}^{\frac{1}{4}}[0,T]}^{2}+\left\Vert \gamma(s^n(t),t)(s^n)'(t)\right\Vert_{W_{2}^{\frac{1}{4}}[0,T]}^{2}+\left\Vert \chi(s^n(t),t)\right\Vert_{W_{2}^{\frac{1}{4}}[0,T]}^{2}\nonumber+\\
+\left\Vert f\right\Vert_{L_{2}(D)}^{2} +  \sum_{k=1}^{n-1}{\bf 1}_{+}(s_{k+1}-s_{k})  \int_{s_{k}}^{s_{k+1}} u^2(x;k) dx \bigg] ,\label{Eq:W:3:36} 
\end{gather}
\end{lemma}

\textbf{Proof}:
By choosing $\eta(x)=2\tau \tilde{u}_{\overline{t}}(x;k)$ in (\ref{Eq:W:1:13}) and by using the following identity
\begin{gather}
2\tau a_{k}(x)\frac{d \tilde{u}(x;k)}{d x}\left( \frac{d \tilde{u}(x;k)}{d x}\right)_{\bar{t}}=a_{k}(x)\left( \frac{d \tilde{u}(x;k)}{d x}\right)^{2}-a_{k-1}(x)\left( \frac{d \tilde{u}(x;k-1)}{d x}\right)^{2}\nonumber\\
-\tau a_{k\bar{t}}(x)\left( \frac{d \tilde{u}(x;k-1)}{d x}\right)^{2}+\tau^{2}a_{k}(x)\left[  \left(\frac{d \tilde{u}(x;k)}{dx}\right)_{\bar{t}}\right]^{2}, \label{Eq:W:3:37}
\end{gather}
we have
\begin{gather}
\int_{0}^{s_{k}}a_{k}(x)\left( \frac{d \tilde{u}(x;k)}{d x}\right)^{2}\, dx-\int_{0}^{s_{k-1}}a_{k-1}(x)\left(\frac{d \tilde{u}(x;k-1)}{d x} \right)^{2}\, dx+2\tau\int_{0}^{s_{k}}(\tilde{u}_{\bar{t}}(x;k))^{2}\, dx+\nonumber\\
+\tau^{2}\int_{0}^{s_{k}}a_{k}(x)\left[\left( \frac{d \tilde{u}(x;k)}{d x}\right)_{\bar{t}} \right]^{2}\leq \tau \int_{0}^{s_{k}}a_{k\bar{t}}(x)\left(\frac{d \tilde{u}(x;k-1)}{d x} \right)^{2}\, dx \nonumber\\ + 2\tau \int_{0}^{s_{k}}b_{k}(x)\frac{d \tilde{u}(x;k)}{d x}\tilde{u}_{\bar{t}}(x;k)\, dx 
+2\tau\int_{0}^{s_{k}}c_{k}(x)\tilde{u}(x;k)\tilde{u}_{\bar{t}}(x;k)\, dx \nonumber\\ -2\tau\int_{0}^{s_{k}}f_{k}(x)u_{\bar{t}}(x;k)\, dx
-2\tau\left[(\gamma_{s^n} (s^n)')^{k}-\chi_{s^n}^{k} \right]\tilde{u}_{\bar{t}}(s_{k};k)-2\tau g_{k}\tilde{u}_{\bar{t}}(0;k)\label{Eq:W:3:38}
\end{gather}
By adding inequalities (\ref{Eq:W:3:38}) with respect to $k$ from 1 to arbitrary $m \leq n$ we derive
\begin{gather}
\int_{0}^{s_{m}}a_{m}(x)\left( \frac{d \tilde{u}(x;m)}{d x}\right)^{2}\, dx+2\tau\sum_{k=1}^{m}\int_{0}^{s_{k}}\tilde{u}_{\bar{t}}^{2}(x;k)\, dx+\tau^{2}\sum_{k=1}^{m}\int_{0}^{s_{k}}a_{k}(x)\left[ \left( \frac{d \tilde{u}(x;k)}{d x}\right)_{\bar{t}}\right]^{2}\, dx \nonumber\\
\leq \tau \sum_{k=1}^{m}\int_{0}^{s_{k}}a_{k\bar{t}}(x)\left( \frac{d \tilde{u}(x;k-1)}{d x}\right)^{2}\, dx +2\tau\sum_{k=1}^{m}\int_{0}^{s_{k}}b_{k}(x)\frac{d \tilde{u}(x;k)}{d x}\tilde{u}_{\bar{t}}(x;k)\, dx \nonumber\\
+2\tau\sum_{k=1}^{m}\int_{0}^{s_{k}}c_{k}(x)\tilde{u}(x;k)\tilde{u}_{\bar{t}}(x;k)\, dx   
 -2\tau\sum_{k=1}^{m}\int_{0}^{s_{k}}f_{k}(x)\tilde{u}_{\bar{t}}(x;k)\, dx\nonumber\\
+\int_{0}^{s_{0}}a_{0}(x)\left(\frac{d \phi}{d x} \right)^{2}\, dx-2\tau\sum_{k=1}^{m}\left[ (\gamma_{s^n} (s^n)')^{k}-\chi_{s^n}^{k}\right]\tilde{u}_{\bar{t}}(s_{k};k)-2\tau \sum_{k=1}^{m} g_{k}\tilde{u}_{\bar{t}}(0;k)\label{Eq:W:3:39}
\end{gather}
By using (\ref{Eq:W:1:6}),(\ref{conditionon_a}) and by applying Cauchy inequalities with appropriately chosen
$\ep > 0$, from (\ref{Eq:W:3:39}) it follows that
\begin{gather}
a_{0}\int_{0}^{s_{m}}\left| \frac{d \tilde{u}(x;k)}{d x}\right|^{2}\,dx+\tau \sum_{k=1}^{m}\int_{0}^{s_{k}}\tilde{u}_{\bar{t}}^{2}(x;k)\, dx + a_{0}\tau^{2}\sum_{k=1}^{m}\int_{0}^{s_{k}}\left[ \left( \frac{d \tilde{u}(x;k)}{d x}\right)_{\bar{t}}\right]^{2}\, dx\leq \nonumber\\
\leq  C \tau \sum_{k=1}^{m}\bigg[\int_{0}^{s_{k}}u^{2}(x;k)\, dx + \int_{0}^{s_{k}}\left|\frac{d u(x;k)}{d x} \right|^{2}\, dx + \int_{0}^{s_{k}}f_{k}^{2}(x)\, dx \bigg]\nonumber\\ C \int_{0}^{s_{0}}\left| \frac{d \phi}{d x}\right|^{2}\, dx -2\tau \sum_{k=1}^{m}\left[(\gamma_{s^n} (s^n)')^{k}-\chi_{s^n}^{k} \right]\tilde{u}_{\bar{t}}(s_{k};k) - 2\tau\sum_{k=1}^{m}g_{k}\tilde{u}_{\bar{t}}(0;k)\, dx\label{Eq:W:3:40}
\end{gather}
where $C$ is independent of $n$. Note that we replaced $\tilde{u}$ with $u$ in first two integrals on the right-hand side of 
(\ref{Eq:W:3:40}). 
Since $\gamma, \chi \in W_2^{1,1}(D)$ we have $\gamma(s^n(t),t), \chi(s^n(t),t) \in W_2^{\frac{1}{4}}[0,T]$ (\cite{Nikolski, BIN, LSU}) and
\begin{equation}\label{Eq:W:3:41}
\Vert \gamma(s^n(t),t)\Vert_{W_2^{\frac{1}{4}}[0,T]}\leq C \Vert \gamma \Vert_{W_2^{1,1}(D)}, \  \Vert \chi(s^n(t),t)\Vert_{W_2^{\frac{1}{4}}[0,T]}  \leq C \Vert \chi \Vert_{W_2^{1,1}(D)},
\end{equation}
where $C$ is independent of $n$. According to Lemma~\ref{mappings} ${\mathcal P}_n([v]_n) \in V_{R+1}$. By applying Morrey inequality to $(s^n)^{\prime}$ we easily deduce that $\gamma(s^n(t),t)(s^n)^{\prime}(t) \in W_2^{\frac{1}{4}}[0,T]$ and moreover,
\begin{equation}\label{Eq:W:3:42}
\Vert \gamma(s^n(t),t)(s^n)^\prime(t)\Vert_{W_2^{\frac{1}{4}}[0,T]}\leq C_1 \Vert \gamma(s^n(t),t)\Vert_{W_2^{\frac{1}{4}}[0,T]} \Vert s^n\Vert_{W_2^2[0,T]} \leq C \Vert \gamma \Vert_{W_2^{1,1}(D)},
\end{equation}
where $C$ is independent of $n$.

Let $w(x,t)$ be a function in $W_{2}^{2,1}(D)$ such that
\begin{equation}\label{Eq:W:3:43}
w(x,0)=\phi(x) \quad \text{for}~x\in [0,s_0], \ a(0,t)w_{x}(0,t)=g(t),\ \quad \text{for a.e.}~t\in [0,T]
\end{equation}
\begin{equation}\label{Eq:W:3:44}
a(s^n(t),t)w_{x}(s^n(t),t)=\gamma(s^n(t),t)(s^n)^{\prime}(t)-\chi(s^n(t),t) \quad \text{for a.e.}~t\in [0,T]
\end{equation}
and
\begin{gather}
\left\Vert w\right\Vert_{W_{2}^{2,1}(D)} \leq C\Big [ \left\Vert g\right\Vert_{W_{2}^{\frac{1}{4}}[0,T]} +\left\Vert \phi(x)\right\Vert_{W_{2}^{1}[0,s_{0}]}\nonumber\\
+\left\Vert \gamma(s^n(t),t)(s^n)^{\prime}(t)-\chi(s^n(t),t)\right\Vert_{W_{2}^{\frac{1}{4}}[0,T]}\Big ]\label{Eq:W:3:45}
\end{gather}
The existence of $w$ follows from the result on traces of Sobolev functions \cite{BIN,Nikolski}.
For example, $w$ can be constructed as a solution from $W_2^{2,1}(\Omega^n)$ of the heat equation 
in 
\[ \Omega^n=\{0<x<s^n(t), 0<t<T \} \]
under initial-boundary conditions (\ref{Eq:W:3:43}),(\ref{Eq:W:3:44})with subsequent continuation to $W_2^{2,1}(D)$ with norm preservation \cite{Solonnikov1, Solonnikov2}. 

Hence, by replacing in the original problem (\ref{Eq:W:1:1})-(\ref{Eq:W:1:4}) $u$ with $u-w$ we can derive modified (\ref{Eq:W:3:40}) without the last three terms on the right-hand side and with $f$, replaced by
\begin{equation}\label{Eq:W:3:46}
F=f+w_{t}-(a w_x)_x-bw_x-cw\in L_{2}(D).
\end{equation}
By using the stability estimation (\ref{firstenergyestimate3}),  from modified (\ref{Eq:W:3:40}),(\ref{Eq:W:3:45}) and (\ref{Eq:W:3:46}), (\ref{Eq:W:3:36}) follows.  Lemma is proved.

In the next lemma we prove (\ref{secondenergyestimate1}) with $l$ being replaced with $s_k$ on the left-hand side.
\begin{lemma}\label{3.4}
For all sufficiently small $\tau$, discrete state vector $[u([v]_n)]_n$ satisfies the following estimation:
\begin{gather}
\max_{1 \leq k \leq n}\int_{0}^{s_{k}}\left| \frac{d u(x;k)}{d x}\right|^{2}\, dx+\tau \sum_{k=1}^{n}\int_{0}^{s_{k}}u_{\bar{t}}^{2}(x;k)\, dx \leq
 \nonumber\\ C \bigg[ \left\Vert \phi\right\Vert_{W_{2}^{1}[0,l]}^{2}+\left\Vert g\right\Vert_{W_{2}^{\frac{1}{4}}[0,T]}^{2}+\left\Vert \gamma(s^n(t),t)(s^n)'(t)\right\Vert_{W_{2}^{\frac{1}{4}}[0,T]}^{2}+\left\Vert \chi(s^n(t),t)\right\Vert_{W_{2}^{\frac{1}{4}}[0,T]}^{2}\nonumber+\\
+\left\Vert f\right\Vert_{L_{2}(D)}^{2} +  \sum_{k=1}^{n-1}{\bf 1}_{+}(s_{k+1}-s_{k})  \int_{s_{k}}^{s_{k+1}} u^2(x;k) dx \bigg] ,\label{Eq:W:3:47} 
\end{gather}
\end{lemma}

\textbf{Proof}: Obviously, we can equivalently replace $\tilde{u}$ with $u$ in the first term on the left-hand side of (\ref{Eq:W:3:36}). We can do so also in the second term provided $s_{k-1} \ge s_k$ for all
$k=\overline{1,m}$. Hence, we only need to estimate
\[ \int\limits_0^{s_k}u_{\bar{t}}^2(x;k)\, dx, \ s_{k-1}<s_k.  \] 
By using (\ref{Eq:W:2:28}) we have
\begin{gather}
\int\limits_{0}^{s_k}u_{\bar{t}}^2(x;k)dx=\int\limits_{0}^{s_{k-1}}\tilde{u}_{\bar{t}}^2(x;k)dx+\int\limits_{s_{k-1}}^{s_k}u_{\bar{t}}^2(x;k)dx\nonumber\\
\int\limits_{s_{k-1}}^{s_k}\Big | \frac{u(x;k)-u(x;k-1)}{\tau}\Big |^2dx \le 2 \int\limits_{s_{k-1}}^{s_k}\Big | \frac{u(x;k)-u(2s_{k-1}-x;k)}{\tau}\Big |^2dx\nonumber\\ + 2 \int\limits_{s_{k-1}}^{s_k}\Big | \frac{u(2s_{k-1}-x;k)-u(2s_{k-1}-x;k-1)}{\tau}\Big |^2dx \le 2 \int\limits_{s_{k-1}}^{s_k}\Big |\frac{1}{\tau} \int\limits_{2s_{k-1}-x}^{x}\frac{du(y;k)}{dy}\Big |^2dx \nonumber\\
+2\int\limits_{s_{k-1}-(s_k-s_{k-1})}^{s_{k-1}}\tilde{u}_{\bar{t}}^2(x;k)dx \le \frac{2}{\tau^2}  \int\limits_{s_{k-1}}^{s_k} \int\limits_{2s_{k-1}-x}^{x}\Big |\frac{du(y;k)}{dy}\Big |^2 dy 2(x-s_{k-1})dx \nonumber\\ + 2\int\limits_{s_{k-1}-C\tau}^{s_{k-1}}\tilde{u}_{\bar{t}}^2(x;k)dx \le
2 \int\limits_{s_{k-1}-C\tau}^{s_{k}}\Big |\frac{d\tilde{u}(x;k)}{dx}\Big |^2dx + 2\int\limits_{s_{k-1}-C\tau}^{s_{k-1}}\tilde{u}_{\bar{t}}^2(x;k)dx.  \label{Eq:W:3:48}
\end{gather}
Hence, for sufficiently small $\tau$ we have
\begin{gather}
\int\limits_{0}^{s_k}u_{\bar{t}}^2(x;k)dx \le 2 \int\limits_{s_{k-1}-C\tau}^{s_{k}}\Big |\frac{d\tilde{u}(x;k)}{dx}\Big |^2dx + \int\limits_{0}^{s_{k-1}}\tilde{u}_{\bar{t}}^2(x;k)dx \nonumber\\
+2\int\limits_{s_{k-1}-C\tau}^{s_{k-1}}\tilde{u}_{\bar{t}}^2(x;k)dx \le 2 \int\limits_{0}^{s_{k}}\Big |\frac{d\tilde{u}(x;k)}{dx}\Big |^2dx + 3\int\limits_{0}^{s_{k-1}}\tilde{u}_{\bar{t}}^2(x;k)dx. \label{Eq:W:3:49}
\end{gather}
From (\ref{Eq:W:3:36}) and (\ref{Eq:W:3:49}), (\ref{Eq:W:3:47}) follows. Lemma is proved.

It can be easily seen that Theorem~\ref{secondenergyestimate} follows from Lemma~\ref{3.4} and extension   Lemma~\ref{extensionlemma}.

Second energy estimate (\ref{secondenergyestimate1}) allows to strengthen the result of Theorem~\ref{continuity1}. 
\begin{theorem}\label{continuity2}
Let $[v]_n \in V_R^n, n=1,2,...$ be a sequence of discrete controls and the sequence $\{\mathcal{P}_n([v]_n\}$ converges strongly in $W_2^1[0,T] \times L_2[0,T]$ to $v=(s,g)$. Then the sequence $\{\hat{u}^\tau\}$ converges as $\tau \to 0$ weakly in $W_2^{1,1}(\Omega)$ to weak solution $u \in W_2^{1,1}(\Omega)$ of the problem (\ref{Eq:W:1:1})-(\ref{Eq:W:1:4}), i.e. to the solution of the integral identity (\ref{Eq:W:1:9}). Moreover, $u$ satisfies the energy estimate
\begin{equation}\label{W211estimate}
\Vert u \Vert_{W_2^{1,1}(D)}^2 \le C \Big ( \Vert \phi \Vert_{W_2^1(0,s_0)}^2 + \Vert g \Vert_{W_2^{\frac{1}{4}}[0,T]}^2 + \Vert f \Vert_{L_2(D)}^2 + 
\Vert \gamma \Vert_{W_2^{1,1}(D)}^2 + \Vert \chi \Vert_{W_2^{1,1}(D)}^2 \Big )
\end{equation}
\end{theorem}

\textbf{Proof}: The last term on the right-hand side of the second energy estimate (\ref{secondenergyestimate1}) is estimated in Theorem~\ref{continuity1} along (\ref{Eq:W:3:19})-(\ref{Eq:W:3:23}). By using Theorems~\ref{firstenergy} and ~\ref{continuity1}, from (\ref{secondenergyestimate1}),(\ref{Eq:W:3:42}) it follows that the sequence $\{\hat{u}^\tau\}$
satisfies the estimate
\begin{equation}\label{Eq:W:3:51}
\Vert \hat{u}^\tau \Vert_{W_2^{1,1}(D)}^2 \le C \Big ( \Vert \phi \Vert_{W_2^1(0,s_0)}^2 + \Vert g \Vert_{W_2^{\frac{1}{4}}[0,T]}^2 + \Vert f \Vert_{L_2(D)}^2 + 
\Vert \gamma \Vert_{W_2^{1,1}(D)}^2 + \Vert \chi \Vert_{W_2^{1,1}(D)}^2 \Big )
\end{equation}
Hence, $\{\hat{u}^\tau\}$ is weakly precompact in $W_2^{1,1}(D)$. It follows that it is strongly precompact
in $L_2(D)$. Let $u$ be a weak limit point of $\{\hat{u}^\tau\}$ in $W_2^{1,1}(D)$, and therefore a strong limit point in $L_2(D)$. From (\ref{secondenergyestimate1}) it follows that
\[ \Vert \hat{u}^\tau - u^\tau \Vert_{L_2(D)}^2=\frac{1}{3}\tau^3 \sum\limits_{k=1}^{n} \int\limits_0^l u_{\bar{t}}^2(x;k)dx \to 0, \quad \text{as}~\tau \to 0. \]
Therefore, $u$ is a strong limit point of the sequence $\{u^\tau\}$ in $L_2(D)$. By Theorem~\ref{continuity1} whole sequence $\{u^\tau\}$
converges weakly in $W_2^{1,0}(\Omega)$ to the unique weak solution from $V_2^{1,0}(\Omega)$ of the problem (\ref{Eq:W:1:1})-(\ref{Eq:W:1:4}). Hence, $u$ is a weak solution of the problem (\ref{Eq:W:1:1})-(\ref{Eq:W:1:4}) and we conclude that whole sequence $\{\hat{u}^\tau\}$ converges weakly in $W_2^{1,1}(D)$ to $u \in W_2^{1,1}(D)$ which is a weak solution of the problem (\ref{Eq:W:1:1})-(\ref{Eq:W:1:4}) from $W_2^{1,1}(\Omega)$. From the property of weak convergence it follows that
$u$ satisfies (\ref{W211estimate}). Theorem is proved.

In particular, Theorem~\ref{continuity2} implies the following  existence result:
\begin{corollary}\label{W211solution}
For arbitrary $v=(s,g)\in V_R$ there exists a weak solution $u \in W_2^{1,1}(\Omega)$ of the problem (\ref{Eq:W:1:1})-(\ref{Eq:W:1:4}) which satisfy the energy estimate (\ref{W211estimate})
\end{corollary}

\textbf{Remark}: In fact, we proved slightly higher regularity of $u$, and both in Theorem~\ref{continuity2} 
and Corollary~\ref{W211solution} $W_2^{1,1}(D)$-norm on the left-hand side of (\ref{W211estimate}) can be replaced with the norm
\[ {\bf \Vert} u {\bf \Vert} = \max_{0\le t \le T}\Vert u(x,t)\Vert_{W_2^1[0,l]}+\Vert u_t \Vert_{L_2(D)} \]
\subsection{Proof of Theorem~\ref{existence} }\label{proofofexistence}
Let $\{ v_{n}\} \in V_{R}$ be a minimizing sequence
\[ \lim_{n\to\infty} \mathcal{J}(v_n)= \mathcal{J}_* \]
Sequence $v_{n}=(s_{n},g_{n})$ is weakly precompact in $W_{2}^{2}[0,T]\times W_{2}^{1}[0,T]$. Assume that the whole sequence $v_{n}=(s_{n},g_{n})$ converge to some limit function $v=(s,g)\in V_R$ weakly in $W_{2}^{2}[0,T] \times W_{2}^{1}[0,T]$, and hence, strongly in $W_{2}^{1}[0,T] \times L_{2}[0,T]$. Let $u_n=u(x,t;v_n), u=u(x,t;v) \in W_{2}^{1,1}(D)$ are weak solutions of (\ref{Eq:W:1:1})-(\ref{Eq:W:1:4})
in $W_2^{1,1}(\Omega_n)$ and $W_2^{1,1}(\Omega)$ respectively. By Corollary~\ref{W211solution}, both 
satisfy energy estimation (\ref{W211estimate}) with $g_n$ and $g$ on the right-hand side respectively.
Since $v_n \in V_R$, $\Vert u_n \Vert_{W_2^{1,1}(D)}$ is uniformly bounded. Hence, the sequence $\Delta u = u_n-u$  it satisfies 
\begin{equation}
\Vert \Delta u \Vert_{W_{2}^{1,1}(D)} \leq C\label{Eq:W:3:52}
\end{equation}
uniformly with respect to $n$. Accordingly, $\{\Delta u\}$ is weakly precompact in $W_{2}^{1,1}(D)$. Without loss of generality assume that the whole sequence $u_{n}-u$ converges weakly in $W_{2}^{1,1}(D)$ to some function $v \in W_{2}^{1,1}(D)$. Let us subtract integral identities (\ref{Eq:W:1:9}) for $u_n$ and $u$,
by assuming that the fixed test function $\Phi$ belongs to $W_2^{1,1}(D)$. Indeed, otherwise $\Phi$ can be continued to $D$ as an element of $W_2^{1,1}(D)$.
\begin{gather}
\int_{0}^{T}\int_{0}^{s(t)}\bigg\{a \Delta u_{x} \Phi_{x} -b \Delta u_{x} \Phi - c \Delta u \Phi + \Delta u_{t}\Phi_{x} \bigg\}\, dx\, dt+\int_{0}^{T}\left( g_{n}(t)-g(t)\right)\Phi(0,t)\, dt \nonumber\\
+\int_{0}^{T}\left[\gamma(s_{n}(t),t)s_{n}'(t)-\gamma(s(t),t)s'(t) - \chi(s_{n}(t),t) + \chi(s(t),t) \right]\Phi(s(t),t)\, dt\nonumber\\
+\int_{0}^{T}\int_{s(t)}^{s_{n}(t)}\left\{ a (u_{n})_x\Phi_{x} -b (u_{n})_x\Phi - c u_{n}\Phi + (u_{n})_t\Phi + f \Phi\right\}\, dx\, dt \nonumber\\+ \int_{0}^{T}\left[\gamma(s_{n}(t),t)s_{n}'(t) - \chi(s_{n}(t),t)\right]\left[ \Phi(s_{n}(t),t)-\Phi(s(t),t)\right]\, dt=0\label{Eq:W:3:53}
\end{gather}
By using energy estimate (\ref{W211estimate}), and continuity of traces $\gamma(s(t),t), \chi(s(t),t)$ of elements $\gamma, \chi \in W_2^{1,1}(D)$, strongly in $L_2[0,T]$ with respect to $s\in W_2^1[0,T]$, passing to the limit as $n \to +\infty$, from (\ref{Eq:W:3:53}) it follows that the weak limit function $v$ satisfies the integral identity
\begin{equation}\label{Eq:W:3:54}
\int_{0}^{T}\int_{0}^{s(t)}\left\{a v_{x}\Phi_{x}-b v_{x}\Phi + c v \Phi + v_t \Phi \right\}\, dx\, dt=0
\end{equation}
for arbitrary $\Phi \in W_{2}^{1,1}(D)$. Since, any element $\Phi \in W_{2}^{1,1}(\Omega)$ can be continued to $D$ as element of $W_{2}^{1,1}(D)$, (\ref{Eq:W:3:54}) is valid for arbitrary $\Phi \in W_{2}^{1,1}(\Omega)$. Hence, $v$ is a weak solution from $W_{2}^{1,1}(\Omega)$ of the problem (\ref{Eq:W:1:1})-(\ref{Eq:W:1:4}) with $f=g=\gamma=\chi=0$. From (\ref{W211estimate}) and uniqueness
it follows that $v=0$. Thus $u_{n}$ converges to $u$ weakly in $W_{2}^{1,1}(D)$. From Sobolev trace theorem (\cite{BIN,Nikolski}) it follows that
\[
\Vert u_{n}(0,t)-u(0,t)\Vert_{L^{2}[0,T]} \to 0, \
 \Vert u_{n}(s(t),t)-u(s(t),t))\Vert_{L^{2}[0,T]} \to 0 \quad \text{as}~ n \to \infty, \]
\[ \Vert u_{n}(s_n(t),t)-u(s(t),t))\Vert_{L^{2}[0,T]} \le \Vert u_{n}(s_n(t),t)-u_n(s(t),t)\Vert_{L^{2}[0,T]} \]
\[ +\Vert u_{n}(s(t),t)-u(s(t),t))\Vert_{L^{2}[0,T]} \to 0 \quad \text{as}~ n \to \infty. \]
Hence, we have
\[ \mathcal{J}(v)=\lim_{n\to\infty}\mathcal{J}(v_n)=\mathcal{J}_{*}\]
and $v$ is a solution of the Problem I. Theorem is proved.

\textbf{Remark}: By applying first and second energy estimates we proved that functional $\mathcal{J}(v)$
is weakly continuous in $W_2^2[0,T] \times W_2^1[0,T]$. Since $V_R$ is weakly compact existence of the optimal control follows from Weierstrass theorem in weak topology.
\subsection{Proof of Theorem~\ref{convergence}}\label{proofofconvergence}
We split the remainder of the proof into three lemmas.
\begin{lemma}\label{condition(3)}
Let $\mathcal{J}_{*}(\pm \epsilon) = \inf\limits_{V_{R\pm \epsilon}} \mathcal{J}(v)$, $\ep >0$. Then
\begin{equation}
\lim_{\epsilon \to 0} \mathcal{J}_{*}(\epsilon) = \mathcal{J}_{*} = \lim_{\epsilon \to 0}\mathcal{J}_{*}(-\epsilon)\label{Eq:W:3:55}
\end{equation}
\end{lemma}
\textbf{Proof:} Note that for $0 < \epsilon_{1}< \epsilon_{2}$ we have 
\[ \mathcal{J}_{*}(\epsilon_{2}) \leq \mathcal{J}_{*}(\epsilon_{1}) \leq \mathcal{J}_{*} \leq \mathcal{J}_{*}(-\epsilon_{1}) \leq \mathcal{J}_{*}(-\epsilon_{2})\]
Therefore $\lim\limits_{\epsilon \to 0}\mathcal{J}_{*}(\epsilon) \leq \mathcal{J}_{*}$ and $\lim\limits_{\epsilon \to 0}\mathcal{J}_{*}(-\epsilon) \geq \mathcal{J}_{*}$ exist. Let us choose $v_{\epsilon} \in V_{R+\epsilon}$ such that
\begin{equation}\label{Eq:W:3:56}
\lim_{\epsilon \to 0}\big( \mathcal{J}(v_{\epsilon})-\mathcal{J}_{*}(\epsilon)\big)=0
\end{equation}
Since $v_{\epsilon}=(s_{\epsilon},g_{\epsilon})$ is weakly precompact in $W_{2}^{2}[0,T]\times W_{2}^{1}[0,T]$, there exists some subsequence $\epsilon'$ such that
\begin{gather*}
s_{\epsilon'}\to s_{*}~\text{weakly in}~W_{2}^{2}[0,T], \ 
g_{\epsilon'}\to g_{*}~\text{weakly in}~W_{2}^{1}[0,T]~\text{as}~\epsilon'\to 0
\end{gather*}
Moreover, we have $v_*=(s_*,g_*)\in V_R$. 
Since $\mathcal{J}(v)$ is weakly continuous it follows that
\begin{equation}\label{Eq:W:3:57}
\lim_{\epsilon' \to 0}\mathcal{J}(v_{\epsilon'}) = \mathcal{J}(v_*).
\end{equation}
From (\ref{Eq:W:3:56}),(\ref{Eq:W:3:57}) it follows that
\[ \lim_{\epsilon' \to 0}\mathcal{J}_{*}(\epsilon')=\mathcal{J}_{*}\]
which implies the first relation in (\ref{Eq:W:3:55}).

To prove the second relation in (\ref{Eq:W:3:55}), take $\epsilon_{0}>0$ and $\tilde{v}=(\tilde{s},\tilde{g}) \in V_{R-\epsilon_{0}}$. Let $\{\alpha_{k}\}$ be a real sequence with $0 < \alpha_{k} < 1$, $\lim\limits_{k \to +\infty} \alpha_{k}=0$ and set
\[ v_{k}=(s_{k},g_{k})=\alpha_{k} \tilde{v} + (1-\alpha_{k}) v_{*}\]
where $\mathcal{J}(v_{*})=\mathcal{J}_{*}$. We have
$v_{k}\in V_{R-\alpha_{k}\epsilon_{0}}$ and $v_k$ converges to $v_*$ strongly in $W_2^2[0,T]\times W_2^1[0,T]$. Since $\mathcal{J}(v)$ is continuous, $v_k$ is a minimizing sequence: 
\[ \lim_{k\to\infty}\mathcal{J}(v_{k})=\mathcal{J}_{*}\]
For fixed $k$ choose arbitrary $\epsilon$ such that $0<\epsilon<\epsilon_{0}\alpha_{k}$. We obviously have

\[\mathcal{J}_{*}(-\epsilon)\leq \mathcal{J}(v_{k}),\quad 0<\epsilon<\epsilon_{0}\alpha_{k}\]
Passing to the limit as $\epsilon \to 0$ we have 
\[\lim_{\epsilon \to 0}\mathcal{J}_{*}(-\epsilon) \leq \mathcal{J}(v_{k})\]
Now we pass to the limit as $k \to +\infty$ and get
\[ \lim_{\epsilon \to 0}\mathcal{J}_{*}(-\epsilon) \leq \mathcal{J}_{*}\]
Since the opposite inequality is obvious, (\ref{Eq:W:3:55}) follows. Lemma is proved.
\begin{lemma}\label{condition(1)}
For arbitrary $v=(s,g) \in V_{R}$,
\begin{equation}\label{Eq:W:3:58}
\lim\limits_{n\to \infty} \mathcal{I}_n(\mathcal{Q}_n(v))=\mathcal{J}(v)
\end{equation}
\end{lemma}
\textbf{Proof:}
Let $v \in V_{R}$, $u=u(x,t;v)$, $\mathcal{Q}_{n}(v)=[v]_{n}$ and $[u([v]_{n})]_{n}$ be a corresponding discrete state vector. In Theorem~\ref{continuity2} it is proved that the sequence $\{\hat{u}^\tau\}$ converges to $u$ weakly in $W_2^{1,1}(\Omega)$. This implies that the sequences of traces $\{\hat{u}^\tau(0,t)\}$ and $\{\hat{u}^\tau(s(t),t)\}$ converge strongly in $L_2[0,T]$ to corresponding traces $u(0,t)$ and $u(s(t),t)$.
Let us prove that that the sequences of traces $\{u^{\tau}(0,t)\}$ and $\{u^{\tau}(s(t),t)\}$ converge strongly in $L^{2}[0,T]$ to traces $u(0,t)$ and $u(s(t),t)$ respectively. By Sobolev embedding theorem (\cite{BIN,Nikolski}) it is enough to prove that the sequences $\{u^\tau\}$ and $\{\hat{u}^\tau\}$ are equivalent in strong topology of $W_2^{1,0}(\Omega)$. In Theorem~\ref{continuity2} it is proved that they are equivalent in strong topology of $L_2(D)$. It remains only to demonstrate that the sequences of 
derivatives $\{u^\tau_x\}$ and $\{\hat{u}^\tau_x\}$ are equivalent in strong topology of $L_2(\Omega)$. 
We have
\begin{equation}\label{Eq:W:3:59}
\Vert u^\tau_x-\hat{u}^\tau_x\Vert_{L_2(\Omega)}^2\le \frac{1}{3}\sum\limits_{k=1}^n \tau^3 \int\limits_0^{\min(s_{k-1}; s_k)} \Big ( \frac{d \tilde{u}(x;k)}{dx} \Big )_{\bar{t}}^2 dx + \Vert u^\tau_x-\hat{u}^\tau_x\Vert_{L_2(\Gamma_n)}^2,
\end{equation}
where $s_k=s^n(t_k)$, $s^n$ is the first component of $\mathcal{P}_n([v]_n)$ and 
\[ \Gamma_n= \cup_{k=1}^{n} \{ t_{k-1}<t\le t_k, \ \min(s_{k-1}; s_k) < x < s(t) \} \]
Since $s^n$ converges to $s$ uniformly on $[0,T]$, it follows that the Lebesgue measure of $\Gamma_n$ converges to zero as $n \to +\infty$. By Theorems~\ref{continuity1} and ~\ref{continuity2} the integrand is uniformly bounded in $L_2(D)$. Therefore, the second term on the right-hand side of (\ref{Eq:W:3:59})
converges to zero as $n \to +\infty$. First term on the right-hand side of (\ref{Eq:W:3:59}) converges to zero due to stability estimation (\ref{Eq:W:3:36}) and the claim is proved.

Let $\nu^{\tau}(t)=\nu^{k}, \ \mu^{\tau}(t)=\mu^{k}$ , if $t_{k-1}< t \leq t_{k}$, $k=1,\ldots,n$. We have
\begin{equation}
\Vert \nu^{k}-\nu \Vert_{L^{2}[0,T]}\to 0, \ \Vert \mu^{k}-\mu \Vert_{L^{2}[0,T]}\to 0~\text{as}~\tau \to 0\label{Eq:W:3:60}
\end{equation}
We estimate the first term in $\mathcal{I}_n(\mathcal{Q}_n(v))$ as follows
\begin{equation}
\beta_{0}\tau \sum_{j=1}^{n}|u(0;k)-\nu^{k}|^{2}=\beta_{0}\sum_{k=1}^{n}\int_{t_{k-1}}^{t_{k}}|u(0;k)-\nu^{k}|^{2}\, dt=\beta_{0}\int_{0}^{T}|u^{\tau}(0,t)-\nu^{\tau}(t)|^{2}\, dt\label{Eq:W:3:61}
\end{equation}
From (\ref{Eq:W:3:60}) it follows that
\begin{equation}
\lim_{n\to\infty}\beta_{0}\tau\sum_{k=1}^{n}|u(0;k)-\nu^{k}|^{2}=\beta_{0}\Vert u(0,t)-\nu(t) \Vert_{L^{2}[0,T]}^{2}\label{Eq:W:3:62}
\end{equation}
We estimate the second term in $\mathcal{I}_n(\mathcal{Q}_n(v))$ as follows
\begin{gather}
\beta_{1}\tau\sum_{k=1}^{n}|u(s_{k};k)-\mu^{k}|^{2}=2\beta_{1}\sum_{k=1}^{n}\int_{t_{k-1}}^{t_{k}}\int_{s(t)}^{s_{k}}\frac{\partial u^{\tau}}{\partial x}\left(u^{\tau}(s(t),t)-\mu^{\tau}(t) \right)\, dx\, dt\nonumber\\ +\beta_{1}\sum_{k=1}^{n}\int_{t_{k-1}}^{t_{k}}|u(s(t);k)-\mu^{k}|^{2}\, dt  +\beta_{1} \sum_{k=1}^{n}\int_{t_{k-1}}^{t_{k}}\left( \int_{s(t)}^{s_{k}}\frac{\partial u^{\tau}}{\partial x}\, dx \right)^{2}\, dt=I_1+I_2+I_3\label{Eq:W:3:63}
\end{gather}
We have
\begin{equation}
\lim\limits_{n\to \infty} I_2=\beta_{1}\int_{0}^{T}|u^{\tau}(s(t),t)-\mu^{\tau}(t)|^{2}\, dt = \beta_{1}\int_{0}^{T}|u(s(t),t)-\mu(t)|^{2}\, dt\label{Eq:W:3:64}
\end{equation}
Since $\Vert (u^\tau)_x\Vert_{L_2(D)}$ and $\Vert u^\tau(s(t),t)-\mu^\tau \Vert_{L_2[0,T]}$ are uniformly bounded, and $\{s^n\}$ converges to $s$ uniformly on $[0,T]$, by applying CBS inequality it easily follows that
\begin{equation}
\lim\limits_{n\to \infty} I_1 = 0, \ \lim\limits_{n\to \infty} I_3 = 0\label{Eq:W:3:65}
\end{equation}
From (\ref{Eq:W:3:63})-(\ref{Eq:W:3:65}) it follows that
\begin{equation}
\lim_{\tau\to 0}\beta_{1}\tau\sum_{k=1}^{n}\left| u(s_{k};k)-\mu^{k}\right|^{2}=\beta_{1}\int_{0}^{T}\left| u(s(t),t)-\mu(t)\right|^{2}\, dt\label{Eq:W:3:66}
\end{equation}
Therefore, from \eqref{Eq:W:3:61} and \eqref{Eq:W:3:66}, \eqref{Eq:W:3:58} follows. Lemma is proved.

\begin{lemma}\label{condition(2)}
For arbitrary $[v]_{n}\in V_{R}^{n}$
\begin{equation}
\lim\limits_{n\to \infty} \Big ( \mathcal{J}(\mathcal{P}_n([v]_n))-\mathcal{I}_n([v]_n)  \Big )= 0\label{Eq:W:3:67}
\end{equation}
\end{lemma}
\textbf{Proof:} Let $[v]_{n}\in V_{R}^{n}$ and $v^{n}=(s^{n},g^{n})=\mathcal{P}_n([v]_n)$. From  Lemma~\ref{mappings}
it follows that the sequence $\{\mathcal{P}_n([v]_n\}$ is weakly precompact in $W_2^2[0,T]\times W_2^1[0,T]$. Assume that the whole sequence converges to $\tilde{v}=(\tilde{s},\tilde{g})$ weakly in $W_2^2[0,T]\times W_2^1[0,T]$. This implies the strong convegence in $W_2^1[0,T] \times L_2[0,T]$. From the well-known property of weak convergence it follows that $\tilde{v} \in V_{R}$. In particular $s^{n}$ converges to $\tilde{s}$ uniformly on $[0,T]$ and we have
\begin{equation}
\lim_{n\to\infty}\max_{0 \leq i \leq n}\left| s^{n}(t_{i})-\tilde{s}(t_{i})\right|=0\label{Eq:W:3:68}
\end{equation}
Let $\mathcal{Q}_{n}(\tilde{v})=[\tilde{v}]_{n}$ We have
\begin{equation}\label{Eq:W:3:69}
\mathcal{I}_{n}\big( [v]_{n}\big)-\mathcal{J}(v^{n}) = \mathcal{I}_{n}\big( [v]_{n}\big)-\mathcal{I}_{n}\big( [\tilde{v}]_{n}\big)+\mathcal{I}_{n}\big( [\tilde{v}]_{n}\big)-\mathcal{J}(\tilde{v}) + \mathcal{J}(\tilde{v})-\mathcal{J}(v^{n})
\end{equation}
In Section~\ref{proofofexistence} we proved the weak continuity of the functional $\mathcal{J}(v)$, i.e. 
\[
\lim_{n\to\infty}\left( \mathcal{J}(\tilde{v})-\mathcal{J}(v^{n})\right)=0. \]
From Lemma~\ref{condition(1)} it follows that
\[
\lim_{n\to\infty}\left( \mathcal{I}_{n}\big( [\tilde{v}]_{n}\big)-\mathcal{J}(\tilde{v})\right) =0. \]
Hence, we only need to prove that
\begin{equation}\label{Eq:W:3:70}
\lim_{n\to\infty}\left( \mathcal{I}_{n}\big( [v]_{n}\big) - \mathcal{I}_{n}\big( [\tilde{v}]_{n}\big)\right)=0
\end{equation}
Let 
\[ [u([v]_n)]_n=(u_n(x;0),u_n(x;1),...,u_n(x;n)), \  [u([\tilde{v}]_n)]_n=(\tilde{u}(x;0),\tilde{u}(x;1),...,\tilde{u}(x;n))  \] 
are corresponding discrete state vectors according to Definition~\ref{discretestatevector}. Let $s_{k}=s^{n}(t_{k})$, $\tilde{s}_{k}=\tilde{s}(t_{k})$ and
\[ \Delta u(x;k)=u_{n}(x;k)-\tilde{u}(x;k)\]
We have
\begin{gather}
\mathcal{I}_{n}([v]_{n})-\mathcal{I}_{n}([\tilde{v}]_{n})=\beta_{0}\sum_{k=1}^{n}\tau \left(u_{n}(0;k)-f_{0}^{k}\right)^{2}+\beta_{1}\sum_{k=1}^{n}\tau \left(u_{n}(s_{k};k)-f_{1}^{k}\right)^{2}\nonumber\\-\beta_{0}\sum_{k=1}^{n}\tau \left(\tilde{u}(0;k)-f_{0}^{k}\right)^{2}-\beta_{1}\sum_{k=1}^{n}\tau \left(\tilde{u}(\tilde{s}_k;k)-f_{1}^{k}\right)^{2}=\beta_{0}\sum_{k=1}^{n}\tau \left( \Delta u(0;k)\right)^{2}+\nonumber\\
+2\beta_{0}\sum_{k=1}^{n}\tau\Delta u(0;k)\left( \tilde{u}(0;k)-f_{0}^{k}\right)+\beta_{1}\sum_{k=1}^{n}\tau\left( \Delta u(s_{k};k) + \tilde{u}(s_{k};k)-\tilde{u}(\tilde{s}_{k};k)\right)^{2}+\nonumber\\
+2\beta_{1}\sum_{k=1}^{n}\tau \left( \Delta u(s_{k};k)+\tilde{u}(s_{k};k) - \tilde{u}(\tilde{s}_{k};k)\right)\left( \tilde{u}(\tilde{s}_{k};k)-f_{1}^{k}\right)\label{Eq:W:3:71}
\end{gather}
From the estimations of Sections~\ref{firstenergyestimate} and ~\ref{secondenergy} it follows that the sequences $\{u_n(x;k)\}$, $\{\tilde{u}(x;k)\}$ are uniformly bounded in $W_2^1[0,l]$. From \eqref{Eq:W:3:68}
it follows that 
\begin{gather}
\sum_{k=1}^{n}\tau\left( \tilde{u}(s_{k};k)-\tilde{u}(\tilde{s}_{k};k)\right)^{2}=\sum_{k=1}^{n}\tau \left| \int_{s_{k}}^{\tilde{s}_{k}}\frac{d \tilde{u}(x;k)}{d x}dx\right|^{2}\nonumber\\
\leq \max_{1 \leq k \leq n}|s_{k}-\tilde{s}_{k}| \sum_{k=1}^{n}\tau\left| \int_{s_{k}}^{\tilde{s}_{k}}\left| \frac{\partial \tilde{u}(x;k)}{\partial x}\right|^{2}\,dx\right| \to 0, \quad \text{as}~ n\to +\infty.\label{Eq:W:3:72}
\end{gather}
From \eqref{Eq:W:3:71} and \eqref{Eq:W:3:72} it follows that in order to prove \eqref{Eq:W:3:70} it is enough to prove that
\begin{equation}
R=\sum_{k=1}^{n}\tau\left[ \big(\Delta u(0;k)\big)^{2} + \big( \Delta u(s_{k};k)\big)^{2}\right] \to 0 ~\text{as}~n\to +\infty\label{Eq:W:3:73}
\end{equation}
By the Morrey inequality we have
\begin{equation}
R\leq C \sum_{k=1}^{n}\tau \left[ \int_{0}^{s_{k}}\big| \Delta u(x;k)\big|^{2}\, dx + \int_{0}^{s_{k}}\left| \frac{d \Delta u(x;k)}{dx}\right|^{2}\, dx\right]\label{Eq:W:3:74}
\end{equation}
where $C$ is independent of $n$. Let us subtract integral identities \eqref{Eq:W:1:13} for $u_n(x;k)$ and $\tilde{u}(x;k)$, by assuming that the fixed test function $\eta$ belongs to $W_2^1[0,l]$. Indeed, otherwise $\eta$ can be continued to $[0,l]$ as a element of $W_2^1[0,l]$:
\begin{gather}
\int_{0}^{s_{k}}\left (a_{k}(x)\frac{d \Delta u}{d x}\frac{d \eta}{d x} - b_{k}(x)\frac{d\Delta u}{d x}\eta(x) - c_{k}(x)\Delta u \eta + \Delta u_{\bar{t}} \eta\right)\, dx -\chi_{\tilde{s}}^{k}\left[\eta(s_{k})-\eta(\tilde{s}_{k}) \right]+\nonumber\\
+\int_{s_{k}}^{\tilde{s}_{k}}\left( a_{k}(x)\frac{d\tilde{u}}{d x}\frac{d \eta}{dx} - b_{k}(x) \frac{d\tilde{u}}{dx}\eta - c_{k}(x)\tilde{u}\eta + f_{k}(x)\eta  + \tilde{u}_{\bar{t}}\eta \right)\, dx+\nonumber\\
+\left[ (\gamma_{s^n} (s^{n})')^{k}-(\gamma_{\tilde{s}} \tilde{s}')^{k}\right]\eta(s_{k})+(\gamma_{\tilde{s}} \tilde{s}')^{k}\left[ \eta(s_{k})-\eta(\tilde{s}_{k})\right] - \left[ \chi_{s^{n}}^{k}-\chi_{\tilde{s}}^{k}\right]\eta(s_{k})=0\label{Eq:W:3:75}
\end{gather}
Our goal now is to derive from \eqref{Eq:W:3:75} that the right-hand side of \eqref{Eq:W:3:74}  converges to zero as $n\to +\infty$. The proof goes along the same lines as the derivation of the first energy estimate in Lemma~\ref{firstenergy-step1}. By choosing $\eta(x)=2\tau\Delta u(x;k)$ in \eqref{Eq:W:3:75}, and by 
using (\ref{Eq:W:1:6}), Cauchy inequalities with appropriately chosen  $\ep >0$, and Morrey inequality \eqref{Morrey} we derive similar to \eqref{Eq:W:3:5}:
\begin{gather}
a_0 \tau \int_{0}^{s_{k}}\left| \frac{d\Delta u(x;k)}{d x}\right|^{2}\,dx+\int_{0}^{s_{k}}\Delta u^{2}(x;k)\, dx -\int_{0}^{s_{k}}\Delta u^{2}(s;k-1)\, dx+ \tau^2 \int_{0}^{s_k}\Delta u _{\bar{t}}^2(x;k)dx  
\nonumber\\ \leq  C_1 \tau \int_{0}^{s_{k}}\Delta u^{2}(x;k)\, dx
 - 2\tau \int_{s_{k}}^{\tilde{s}_{k}}\bigg[a_{k}(x) \frac{d \tilde{u}}{dx} \frac{d \Delta u(x;k)}{dx} -b_{k}(x)\frac{d \tilde{u}}{dx}\Delta u(x;k) -c_{k}(x)\tilde{u}\Delta u(x;k) \nonumber\\
+f_{k}(x)\Delta u(x;k)+\tilde{u}_{\bar{t}}\Delta u(x;k)\bigg]\, dx-2\int_{t_{k-1}}^{t_{k}}\left(\gamma(s^{n}(t),t)(s^{n})'(t) - \gamma(\tilde{s}(t),t)\tilde{s}'(t) \right)dt \ \Delta u(s_{k};k)\nonumber\\
-2\int_{t_{k-1}}^{t_{k}}\int_{\tilde{s}_{k}}^{s_{k}}\gamma(\tilde{s}(t),t)\tilde{s}'(t)\frac{d \Delta u(x;k)}{dx}\, dx\, dt+2\int_{t_{k-1}}^{t_{k}}\left(\chi(s^{n}(t),t)-\chi(\tilde{s}(t),t) \right)\, dt\Delta u(s_{k};k) \nonumber\\+ 2\int_{t_{k-1}}^{t_{k}}\int_{\tilde{s}_{k}}^{s_{k}}\chi(\tilde{s}(t),t)\frac{d\Delta u(x;k)}{dx}\, dx\, dt\label{Eq:W:3:76}
\end{gather}
By applying the technique along with \eqref{Eq:W:3:6}-\eqref{Eq:W:3:10} from \eqref{Eq:W:3:76} it follows that for all sufficiently small $\tau$
\begin{equation}
\max_{1 \leq k \leq n}\int_{0}^{s_{k}}\Delta u^{2}(x;k)\, dx \leq C \Big ( \sum_{j=1}^{n-1} {\bf 1}_{+}(s_{j+1}-s_{j})\int_{s_{j}}^{s_{j+1}} \Delta u^2(x;j) dx + \sum_{j=1}^{n}|\mathcal{L}_j|\Big )\label{Eq:W:3:77}
\end{equation}
where $C$ is independent of $\tau$ and
\begin{gather}
\mathcal{L}_j=\tau \int_{s_{j}}^{\tilde{s}_{j}}\Big( a_{j}(x)\frac{d \tilde{u}}{dx} \frac{d \Delta u(x;j)}{dx} - b_{j}(x)\frac{d \tilde{u}}{dx} \Delta u(x;j) - c_{j}(x)\tilde{u}\Delta u(x;j) + f_{j}(x)\Delta u(x;j)\nonumber\\
+ \tilde{u}_{\bar{t}}\Delta u(x;j)\Big)\, dx+\int_{t_{j-1}}^{t_{j}}\left[ \gamma(s^{n}(t),t)(s^{n})'(t)-\gamma(\tilde{s}(t),t)\tilde{s}'(t)\right]dt \ \Delta u(s_{j};j) \nonumber\\
+\int_{t_{j-1}}^{t_{j}}\int_{\tilde{s}_{j}}^{s_{j}}\gamma(\tilde{s}(t),t)\tilde{s}'(t)\frac{d \Delta u(x;j)}{dx}\, dx\, dt+\int_{t_{j-1}}^{t_{j}}\left(\chi(s^{n}(t),t)-\chi(\tilde{s}(t),t) \right)dt \  \Delta u(s_{j};j)\nonumber\\
+\int_{t_{j-1}}^{t_{j}}\int_{\tilde{s}_{j}}^{s_{j}}\chi(\tilde{s}(t),t)\frac{d \Delta u(x;j)}{d x}\, dx \, dt\label{Eq:W:3:78}
\end{gather}
Having \eqref{Eq:W:3:77} we perform summation in \eqref{Eq:W:3:76} with respect to $k$ from $1$ to $n$ 
and derive
\begin{gather}
\max_{1 \leq k \leq n}\int_{0}^{s_{k}}\Delta u^{2}(x;k)\, dx + \sum_{k=1}^{n}\tau\int_{0}^{s_{k}}\left| \frac{d \Delta u(x;k)}{d x}\right|^{2}\,dx+\sum_{k=1}^{n}\tau^{2}\int_{0}^{s_{k}}\Delta u_{\bar{t}}^{2}(x;k)\, dx\nonumber\\ \leq C_1 \Big ( \sum_{j=1}^{n-1} {\bf 1}_{+}(s_{j+1}-s_{j})\int_{s_{j}}^{s_{j+1}} \Delta u^2(x;j) dx + \sum_{j=1}^{n}|\mathcal{L}_j|\Big )\label{Eq:W:3:79}
\end{gather}
where $C_1$ is independent of $\tau$. Our next goal is to absorb the first term on the right-hand side of \eqref{Eq:W:3:79} into the left-hand side. We apply the same method used in the proof of Theorem~\ref{continuity1} (see \eqref{Eq:W:3:19}-\eqref{Eq:W:3:25}). The only difference is that in the estimations \eqref{Eq:W:3:22} and \eqref{Eq:W:3:23} we replace $D$ with
\[ \Omega^1_n=\{0<t<T, 0<x<\tilde{s}_1^n(t) \}  \]
Let us also introduce the region
\[ \Omega_n=\bigcup_{k=1}^n \{t_{k-1}<t\leq t_k, 0<x<s_k\} \]
Note that 
\[ \Vert \Delta u^\tau \Vert_{V_2^{1,0}(\Omega_n)}^2 = \max_{1 \leq k \leq n}\int_{0}^{s_{k}}\Delta u^{2}(x;k)\, dx + \sum_{k=1}^{n}\tau\int_{0}^{s_{k}}\left| \frac{d \Delta u(x;k)}{d x}\right|^{2}\,dx \]
Hence, by applying the method used in Theorem~\ref{continuity1} we derive from \eqref{Eq:W:3:79} the following estimate:
\begin{equation}\label{Eq:W:3:80}
\Vert \Delta u^\tau \Vert_{V_2^{1,0}(\Omega_n)}^2 \leq C_2 \Big ( \Vert \Delta u^\tau \Vert_{V_2^{1,0}(\Omega^1_n-\Omega_n)}^2 + \sum_{j=1}^{n}|\mathcal{L}_j| \Big )
\end{equation}
where $C_2$ is independent of $\tau$. Since the sequences $\{\tilde{s}_1^n\}$ and $\{s^n\}$ are equivalent in strong topology of $W_2^1[0,T]$, the first term on the right-hand side of \eqref{Eq:W:3:80} converges to zero as $n\to+\infty$. It only remains to prove that 
\begin{equation}\label{Eq:W:3:81}
\lim\limits_{n\to +\infty}\sum_{j=1}^{n}|\mathcal{L}_j| = 0. 
\end{equation}
We have
\begin{gather}
\sum_{j=1}^{n}\left|\tau \int_{s_{j}}^{\tilde{s}_{j}}a_{j}(x)\frac{d \tilde{u}(x;j)}{d x} \frac{d \Delta u(x;j)}{d x}\, dx\right| = \sum_{j=1}^{n}\left|\int_{t_{j-1}}^{t_{j}}\int_{s_{j}}^{\tilde{s}_{j}}a(x,t)\frac{d \tilde{u}^{\tau}}{d x} \frac{d \Delta u^{\tau}}{d x}\, dx\right|\nonumber\\
\leq C \left\Vert \frac{\partial \tilde{u}^{\tau}}{\partial x}\right\Vert_{L_{2}(\tilde{\Delta})}\left\Vert \frac{\partial \Delta u^{\tau}}{\partial x}\right\Vert_{L_{2}(\tilde{\Delta})}\label{Eq:W:3:82}
\end{gather}
where
\[ \tilde{\Delta}=\bigcup_{j=1}^{n}\{ (x,t):t_{j-1}<t<t_{j},~\min(s_{j},\tilde{s}_{j})<x<\max(s_{j},\tilde{s}_{j})\}\]
From \eqref{Eq:W:3:68} it follows that the Lebesgue measure of $\tilde{\Delta}$ converges to zero as $n\to \infty$. Since by the first energy estimate $W_{2}^{1,0}(D)$ norm of $\tilde{u}^{\tau}$ and $\Delta u^{\tau}$ are uniformly bounded, the right-hand side of \eqref{Eq:W:3:82} converges to zero as $n\to\infty$. For the same reason, the next three terms in the expression of $\sum_{j=1}^{n}|\mathcal{L}_j|$ also converge to zero as $n\to \infty$. We have
\begin{gather}
\sum_{j=1}^{n}\left|\tau\int_{s_{j}}^{\tilde{s}_{j}} \tilde{u}_{\bar{t}}(x;j) \Delta u(x;j)\, dx\right|
\leq  \left( \sum_{j=1}^{n}\tau \int_{0}^{l}\tilde{u}_{\bar{t}}^{2}(x;j)\, dx\right)^{\frac{1}{2}} \Vert\Delta u^{\tau}\Vert_{L_{2}(\tilde{\Delta})}\label{Eq:W:3:83}
\end{gather}
By the second energy estimate $W_2^{1,1}(D)$ norm of $\hat{\tilde{u}}^\tau$ is uniformly bounded. Accordingly, the right-hand side of \eqref{Eq:W:3:83} converges to zero as $n\to +\infty$. 
We have
\begin{gather}
\sum_{j=1}^{n}\int_{t_{j-1}}^{t_{j}}\Big| \Big (\gamma(s^{n}(t),t)(s^{n})'(t)-\gamma(\tilde{s}(t),t)\tilde{s}'(t)\Big ) \Delta u(s_{j},j)\Big |\, dt \leq \nonumber\\
\sum_{j=1}^{n}\int_{t_{j-1}}^{t_{j}}\Big | \Big ( \gamma(s^{n}(t),t)(s^{n})'(t)-\gamma(\tilde{s}(t),t)\tilde{s}'(t)\Big ) \Big ( \Delta u(\tilde{s}(t);j)+\int_{\tilde{s}(t)}^{s_{j}}\frac{\partial \Delta u^{\tau}(x,t)}{\partial x}dx \Big )\Big |dt\nonumber\\
\leq \Vert\gamma(s^{n}(t),t)-\gamma(\tilde{s}(t),t)\Vert_{L_{2}[0,T]}\Vert (s^{n})'\Vert_{C[0,T]}\Vert \Delta u^{\tau}(\tilde{s}(t),t)\Vert_{L_{2}[0,T]} \nonumber\\ 
 +\Vert (s^{n})'(t)-\tilde{s}'(t)\Vert_{L_{2}[0,T]} \Vert\gamma(\tilde{s}(t),t)\Vert_{L_{4}[0,T]}\Vert\Delta u^{\tau}(\tilde{s}(t),t)\Vert_{L_{4}[0,T]}\nonumber\\
+ \sum_{j=1}^{n}\left( \int_{t_{j-1}}^{t_{j}}|\gamma(s^{n}(t),t)(s^{n})'(t)-\gamma(\tilde{s}(t),t)\tilde{s}'(t)|^{2}\, dt\right)^{\frac{1}{2}}\left( \int_{t_{j-1}}^{t_{j}}\left|\int_{\tilde{s}(t)}^{s_{j}}\frac{\partial \Delta u^{\tau}}{\partial x} \, dx\right|^{2}\right)^{\frac{1}{2}}\nonumber\\
\leq C \Vert \gamma(s^{n}(t),t)-\gamma(\tilde{s}(t),t) \Vert_{L_{2}[0,T]} \Vert (s^{n})' \Vert_{W_{2}^{1}[0,T]} \Vert \Delta u^{\tau} \Vert_{W_{2}^{1,0}(D)}\nonumber\\
+ C \Vert (s^{n})'(t)-\tilde{s}'(t) \Vert_{L_{2}[0,T]} \Vert \gamma \Vert_{W_{2}^{1,1}(D)} \Vert\Delta u^{\tau} \Vert_{V_{2}(D)}\nonumber\\ 
+\Vert\gamma(s^{n}(t),t)(s^{n})'(t)-\gamma(\tilde{s}(t),t)\tilde{s}'(t) \Vert_{L^{2}[0,T]} 
 \left\Vert \frac{\partial \Delta u^{\tau}}{\partial x}\right\Vert_{L_{2}(\tilde{\Delta})} \Big (\max_{0 \leq t \leq T}|\tilde{s}(t)-s^{n}(t)|\Big )^{\frac{1}{2}}\label{Eq:W:3:84}
\end{gather}
and all three terms on the right-hand side converge to zero as $n \to +\infty$.
Similarly one can prove that all the last three terms in the expression of $\sum_{j=1}^{n}|\mathcal{L}_j|$ converges to zero as $n\to \infty$. Hence, \eqref{Eq:W:3:81} is proved. 
From \eqref{Eq:W:3:80} and \eqref{Eq:W:3:81}, \eqref{Eq:W:3:73} follows. Lemma is proved. 

Having Lemmas~\ref{condition(3)}, ~\ref{condition(1)} and ~\ref{condition(2)}, Theorem~\ref{convergence} 
follows from Lemma~\ref{generalcriteria}



\end{document}